\pdfoutput=1
\documentclass[11pt,a4paper]{article}
\usepackage[T1]{fontenc}
\usepackage[utf8]{inputenc}
\usepackage{amsmath,amssymb,amsthm}
\usepackage{txfonts}
\usepackage{mathrsfs,mathtools,array,booktabs,enumitem,microtype,needspace,verbatim}
\usepackage[margin=1in]{geometry}
\usepackage[colorlinks=true,linkcolor=black,citecolor=black,urlcolor=blue]{hyperref}
\hypersetup{
 pdftitle={Positive semidefiniteness of h-Shuhan matrices (analogs of Borcherds Cartan matrices) and associated determinant sequences},
 pdfauthor={Weicai Wu and Mingxuan Yang},
 pdfsubject={Positive semidefiniteness of h-Shuhan matrices and determinant sequences}}
\numberwithin{equation}{section}

\allowdisplaybreaks[2]
\setlist[itemize]{leftmargin=2.2em,itemsep=1.5pt,topsep=3pt,parsep=0pt}
\setlist[enumerate]{leftmargin=2.35em,itemsep=1.5pt,topsep=3pt,parsep=0pt}
\renewcommand{\arraystretch}{1.14}
\newlist{parts}{enumerate}{1}
\setlist[parts]{label=\textup{(\roman*)},ref=(\roman*),leftmargin=2.4em,itemsep=2pt,topsep=3pt}
\newlist{mainparts}{enumerate}{1}
\setlist[mainparts]{label=\textup{(\arabic*)},ref=(\arabic*),leftmargin=2.4em,itemsep=3pt,topsep=4pt}

\newcommand{\psup}[1]{^{\langle #1\rangle}}
\newcommand{\affpsup}[2]{^{(#1)\langle #2\rangle}}

\def\dr{\ar@{->}[r]}

\theoremstyle{plain}
\newtheorem{Theorem}{Theorem}[section]
\newtheorem{Lemma}[Theorem]{Lemma}
\newtheorem{Corollary}[Theorem]{Corollary}
\newtheorem{Proposition}[Theorem]{Proposition}
\newtheorem{Conjecture}[Theorem]{Conjecture}
\theoremstyle{definition}
\newtheorem{Definition}[Theorem]{Definition}
\newtheorem{Question}[Theorem]{Question}
\newtheorem{Example}[Theorem]{Example}
\theoremstyle{remark}
\newtheorem{Remark}[Theorem]{Remark}
\begin{document}

\title{Positive semidefiniteness of $h$-Shuhan matrices
(analogs of Borcherds Cartan matrices) and associated determinant sequences}
\author{Weicai Wu \qquad Mingxuan Yang\\[5pt]
\small School of Mathematics and Statistics, Guangxi Normal University,\\[-1pt]
\small Guilin 541004, Guangxi, People's Republic of China\\[3pt]
\small Correspondence: \texttt{weicaiwu@hnu.edu.cn}}
\date{}
\maketitle

\begin{abstract}
We introduce an $h$-Shuhan matrix as a real square matrix $H=(h_{ij})$ with a common diagonal entry $h\ge0$, nonpositive integral off-diagonal entries, and the following Cartan-type condition on each unordered pair $\{i,j\}$: either $h_{ij}=h_{ji}$, or $\max\{h_{ij},h_{ji}\}=-1$ and $\min\{h_{ij},h_{ji}\}\le -2$. These conditions imply that the scalar shift $H+(2-h)I_n$ is a generalized Cartan matrix; not every generalized Cartan matrix arises from an $h$-Shuhan matrix in this way. We classify, for $0\le h<2$, all indecomposable $h$-Shuhan matrices that are positive semidefinite, generalized positive semidefinite, or virtual positive semidefinite (the last notion means that all principal minors are nonnegative). For $h<2$, virtual positive semidefiniteness forces $H+(2-h)I_n$ to be of finite Cartan type, the positive-type case of Vinberg's trichotomy (1971). At $h=2$, we describe the boundary case in terms of finite and affine Cartan matrices. The associated determinant recurrences hold for general $h$. We study congruences for integral $h$ and asymptotic and Fibonacci-type identities for $h>2$.
\end{abstract}


\section{Introduction}\label{s0}
Cartan matrices are fundamental in the classification of finite-dimensional complex semisimple Lie algebras and in the construction of Kac--Moody algebras; see \cite{Hu72,Ka90}. Positivity questions for matrices with Cartan-type sign patterns also occur in the theory of symmetrizable and quasi-Cartan matrices, $Z$-matrices, and semimonotone matrices; see, for example, \cite{BGZ06,HK03,TW19,MM21,N97,B07}.

The starting point is a deformation of the usual Cartan normalization: the diagonal entry $2$ is replaced by a common real parameter $h\ge0$, while the off-diagonal entries remain integral and nonpositive. The name ``Shuhan'' comes from the diagrams introduced by Wu, Liu, and Xie \cite{WLX22} and extended by Wu and Zhou \cite{WZ25}. Both Chinese papers write the name in Latin letters. We retain their terminology and give the required matrix conditions explicitly in Definition~\ref{1.7}.

The diagrammatic motivation can be stated without assuming familiarity with those papers. In the notation of \cite[Definition~1.2]{WZ25}, each vertex carries a positive integer $x_i$, and the directed edge multiplicities $\lambda_{ij}$ satisfy
\[
  2x_i=\sum_{j\ne i}\lambda_{ij}x_j.
\]
Thus, for~$C=(c_{ij})$ defined by~$c_{ii}:=2$ and~$c_{ij}:=-\lambda_{ij}$ for~$i\ne j$, the vertex labels form a positive vector satisfying $Cx=0$. The connected generalized diagrams are related to affine Dynkin diagrams in this way. We study the associated matrices with diagonal $h$ without requiring such a positive integral null vector. In particular, $h=2$ restores the diagonal normalization used in \cite{WLX22,WZ25}; it does not assert that every matrix in Definition~\ref{1.7} is an affine diagram matrix.

For comparison, Borcherds generalized Kac--Moody theory permits imaginary simple roots and more general diagonal behavior \cite{Bor88}, while Chapovalov--Chapovalov--Lebedev--Leites allow normalized diagonal entries $0,1,2$ for Lie-superalgebra Cartan matrices \cite{CCLL10}. Our definition is different: all diagonal entries equal the same real parameter $h$, and
\[
   C(H):=H+(2-h)I_n
\]
is a generalized Cartan matrix. In particular, an $h$-Shuhan matrix need not itself be a Borcherds--Cartan matrix. If~$0\le h<2$ and~$H$ is virtual positive semidefinite, Proposition~\ref{3.3}(ii) implies that~$C(H)$ is of finite type. Hence, the classification of such matrices reduces to the finite Dynkin types, and the determinants of the corresponding shifted families satisfy the recurrences studied in Sections~\ref{s3} and~\ref{s4}.

Vinberg's theory of linear reflection groups \cite{Vin71} provides the appropriate classical comparison. His Cartan matrices have diagonal $2$, nonpositive off-diagonal entries, and the zero-pattern condition $c_{ij}=0$ if and only if $c_{ji}=0$.
For each pair of distinct indices $i,j$, the product of the two opposite off-diagonal entries satisfies
\[
 c_{ij}c_{ji}\in
 \bigl\{4\cos^2(\tfrac{\pi}{m}):m\in\mathbb Z,\ m\ge2\bigr\}\cup[4,\infty).
\]
The products $c_{ij}c_{ji}$ for~$C(H)=(c_{ij})$ satisfy this condition: they are nonnegative integers, and the values $0,1,2,3$ correspond to $m=2,3,4,6$, respectively. Vinberg's trichotomy into positive, zero, and negative types \cite{Vin71} therefore applies to $C(H)$.
 Section~\ref{s2} spells out its spectral relation to our positivity conditions.

The finite-type reduction and the virtual-semidefinite classification are consequences of this classical theory and the finite Dynkin classification. Our purpose is to determine the thresholds for the fixed integral representatives, including the stronger condition $x^{\mathsf T}Hx\ge0$. This condition depends on the symmetric part $\tfrac12(H+H^{\mathsf T})$ and is not invariant under similarity with respect to a positive diagonal matrix, namely a transformation of the form

\[
 H\longmapsto D^{-1}HD,\qquad
 D=\operatorname{diag}(d_1,\ldots,d_n),\qquad d_i>0.
\]
 The explicit nonsymmetric thresholds in Proposition~\ref{4.8}, the affine bounds in Proposition~\ref{4.11}, and their comparison with the virtual thresholds are the results of matrix computations developed here.

Section~\ref{s4} studies recurrences in $n$, closed formulas, and congruences modulo primes for $a_n\psup{h},b_n\psup{h},d_n\psup{h}$; these are classical second-order recurrence identities \cite{S24}. We then study the two-index determinant sequences $a_{i,n}\psup{h},b_{i,n}\psup{h}$, their invariants, the exceptional limits in Lemma~\ref{5.1}, and the factorizations in Theorem~\ref{5.6}.

Let~$\mathbb N$, $\mathbb Z$, $\mathbb R_{\ge0}$, and~$\mathbb R_{>0}$ denote the positive integers, integers, nonnegative real numbers, and positive real numbers, respectively. Put~$[n]\mathrel{:=}\{1,\ldots,n\}$ and let~$\mathbb S_n$ be the permutation group on $[n]$. The identity matrix of order $n$ is denoted by~$I_n$. For a real $n\times n$ matrix~$H$ and~$J\subseteq[n]$, let~$H[J]$ be the principal submatrix with row and column set $J$; by convention $\det H[\varnothing]\mathrel{:=}1$. When $J=\{i_1<\cdots<i_m\}$, we also write $\widetilde H_{i_1,\ldots,i_m}$ for~$H[J]$. We write
\[
   \widehat H:=\frac{H+H^{\mathsf T}}2.
\]
For a symmetric real matrix, ``positive semidefinite'' and ``positive definite'' have their standard quadratic-form meanings. We write $Q\succeq0$ when the real symmetric matrix~$Q$ is positive semidefinite.

\Needspace{9\baselineskip}
\begin {Definition} \label{1.1}
Let~$H\in M_n(\mathbb R)$.
\begin{itemize}
\item $H$ is \emph{virtual positive semidefinite} if~$\det H[J]\ge0$ for every nonempty $J\subseteq[n]$; it is \emph{virtual positive definite} if all these determinants are strictly positive.
\item $H$ is \emph{generalized positive semidefinite} if~$x^{\mathsf T}Hx\ge0$ for every $x\in\mathbb R^n$; it is \emph{generalized positive definite} if the inequality is strict for every $x\ne0$.
\end{itemize}
\end {Definition}

For clarity, the three positivity notions used in the paper are compared below.
\[
\begin{array}{c|c|c}
\text{notion} & \text{criterion} & \text{symmetry required}\\ \hline
\text{positive semidefinite} & x^{\mathsf T}Hx\ge0 & H=H^{\mathsf T}\\
\text{generalized positive semidefinite} & x^{\mathsf T}Hx\ge0 & \text{no}\\
\text{virtual positive semidefinite} & \det H[J]\ge0\ (\varnothing\ne J\subseteq[n]) & \text{no}
\end{array}
\]
For symmetric $H$, the three semidefinite notions coincide by the principal-minor criterion. For nonsymmetric $H$, generalized positive semidefiniteness implies virtual positive semidefiniteness (Proposition~\ref{1.4}), but the converse fails.

\begin {Lemma} \label{1.2}
A real matrix~$H$ is generalized positive semidefinite if and only if~$\widehat H=\tfrac{1}{2}(H+H^{\mathsf T})$ is positive semidefinite.
\end {Lemma}
\begin{proof}
Write $H=\widehat H+K$, where~$K=\tfrac{1}{2}(H-H^{\mathsf T})$ is skew-symmetric. For every real vector~$x$, the scalar $x^{\mathsf T}Kx$ equals its transpose $-x^{\mathsf T}Kx$, and hence is zero. Therefore, \[
  x^{\mathsf T}Hx=x^{\mathsf T}\widehat Hx
\]
for every $x\in\mathbb R^n$. The two nonnegativity conditions are consequently equivalent. \end{proof}

\begin {Lemma} \label{1.3}
Assume that~$G$ is a real skew-symmetric matrix of order $n$.

{\rm (i)} If~$H$ is positive definite, then $\det(H+G)\ge \det H$, with equality if and only if~$G=0$.

{\rm (ii)} If~$H$ is positive semidefinite, then $\det(H+G)\ge \det H$.
\end {Lemma}
\begin{proof}
{\rm (i)} Since~$H$ is positive definite, there is an invertible $P$ such that~$P^{\mathsf T}HP=I_n$. The matrix~$P^{\mathsf T}GP$ is skew-symmetric, so an orthogonal change of basis puts it into blocks
$\left(\begin{smallmatrix}0&t_j\\-t_j&0\end{smallmatrix}\right)$ and possibly one zero block. Hence, \[
 \frac{\det(H+G)}{\det H}
 =\det\!\left(I_n+P^{\mathsf T}GP\right)
 =\prod_j(1+t_j^2)\ge1.
\]
Equality holds exactly when every $t_j=0$, equivalently $G=0$.

{\rm (ii)} For every $t>0$, the matrix~$H+tI_n$ is positive definite, so part (i) gives

\[\det(H+tI_n+G)\ge\det(H+tI_n).\]
 Letting $t\to0^+$ and using continuity of the determinant proves the claim. \end{proof}

\begin {Proposition} \label{1.4}
If~$H$ is generalized positive semidefinite, then

{\rm (i)} every principal submatrix $H[J]$ is generalized positive semidefinite;

{\rm (ii)} $H$ is virtual positive semidefinite. The converse is false in general.
\end {Proposition}
\begin{proof}
{\rm (i)} Fix $J\subseteq[n]$ and~$y\in\mathbb R^{|J|}$. Embed $y$ into a vector~$x\in\mathbb R^n$ by placing the entries of $y$ in the coordinates indexed by~$J$ and putting zero in the remaining coordinates. Then, $x^{\mathsf T}Hx=y^{\mathsf T}H[J]y\ge0$.

{\rm (ii)} By part (i) and Lemma~\ref{1.2}, the matrix~$\widehat{H[J]}$ is positive semidefinite for every nonempty $J$. Decompose
\[
 H[J]=\widehat{H[J]}+\frac{H[J]-H[J]^{\mathsf T}}2.
\]
The second summand is skew-symmetric, so Lemma~\ref{1.3}(ii) yields
$\det H[J]\ge\det\widehat{H[J]}\ge0$. Thus, every principal minor is nonnegative.

For the converse, take
$H=\left(\begin{smallmatrix}2&-\tfrac72\\-1&2\end{smallmatrix}\right)$. Its principal minors are $2,2,\tfrac12>0$, so it is virtual positive definite, whereas $(1,1)H(1,1)^{\mathsf T}=-\tfrac{1}{2}<0$. \end{proof}

\begin {Remark} \label{1.5}
{\rm (i)} If~$H$ and~$G$ are generalized positive semidefinite, then so is $H+G$. In particular, $H+\lambda I_n$ is generalized positive definite for every $\lambda>0$ whenever $H$ is generalized positive semidefinite.

{\rm (ii)} The analogous assertion for virtual positive semidefiniteness is false without further hypotheses. For example,
$H=\left(\begin{smallmatrix}2&-\tfrac72\\-1&2\end{smallmatrix}\right)$ and
$G=\left(\begin{smallmatrix}2&-1\\-\tfrac72&2\end{smallmatrix}\right)$ are virtual positive definite, but $H+G=\left(\begin{smallmatrix}4&-\tfrac92\\-\tfrac92&4\end{smallmatrix}\right)$ has determinant $-\tfrac{17}{4}<0$ and hence is not virtual positive semidefinite.
\end {Remark}

For each $\sigma\in\mathbb S_n$, set $H_\sigma:=(h_{\sigma(i),\sigma(j)})_{1\le i,j\le n}$.

\begin {Lemma} \label{1.6}
For every $\sigma\in\mathbb S_n$:

{\rm (i)} $\det H=\det H_\sigma$;

{\rm (ii)} $H$ is virtual positive semidefinite if and only if~$H_\sigma$ is virtual positive semidefinite;

{\rm (iii)} $H$ is generalized positive semidefinite if and only if~$H_\sigma$ is generalized positive semidefinite.
\end {Lemma}
\begin{proof}
Let~$P_\sigma$ be the permutation matrix corresponding to $\sigma$. Then, $H_\sigma=P_\sigma^{\mathsf T}HP_\sigma$, so

\[
 \det H_\sigma=(\det P_\sigma)^2\det H=\det H.
\]
This proves (i).

For every $J\subseteq[n]$, the principal submatrix $H_\sigma[J]$ is permutation equivalent to $H[\sigma(J)]$; hence, their determinants are equal. As $J\mapsto\sigma(J)$ permutes the subsets of $[n]$, this proves (ii). Finally, for every $x\in\mathbb R^n$,
\[
 x^{\mathsf T}H_\sigma x=(P_\sigma x)^{\mathsf T}H(P_\sigma x),
\]
and~$P_\sigma$ is bijective on $\mathbb R^n$, which proves (iii). \end{proof}

\begin {Proposition} \label{1.8}
Assume $n\ge2$. Then, $H$ is virtual positive semidefinite if and only if every principal submatrix $H[K]$ with~$|K|=n-1$ is virtual positive semidefinite and~$\det H\ge0$.
\end {Proposition}
\begin{proof}
The forward implication is immediate. Conversely, let~$\varnothing\ne J\subseteq[n]$. If~$J=[n]$, then $\det H[J]=\det H\ge0$. If~$|J|\le n-1$, choose $K\subseteq[n]$ with~$J\subseteq K$ and~$|K|=n-1$. Then, $H[J]$ is a principal submatrix of the virtual positive semidefinite matrix~$H[K]$, so~$\det H[J]\ge0$. Thus, all principal minors of $H$ are nonnegative. \end{proof}

\begin {Definition} \label{1.7}
An $n\times n$ real matrix~$H=(h_{ij})$ is called an \emph{$h$-Shuhan matrix} if

{\rm (S1)} $h_{ii}=h$ for every $i$, where~$h\in\mathbb R_{\ge0}$;

{\rm (S2)} $h_{ij}\in\mathbb Z_{\le0}$ for every $i\ne j$;

{\rm (S3)} for every unordered pair $\{i,j\}$ with~$i\ne j$, either $h_{ij}=h_{ji}$, or
\[
   \max\{h_{ij},h_{ji}\}=-1,
   \qquad \min\{h_{ij},h_{ji}\}\le -2.
\]
Equivalently, in the asymmetric case $\{h_{ij},h_{ji}\}=\{-1,-m\}$ for some integer $m\ge2$.
\end {Definition}

Condition (S3) is deliberately invariant under interchanging $i$ and~$j$. Hence, simultaneous permutation of rows and columns preserves the class: every principal submatrix and every $H_\sigma$ is again an $h$-Shuhan matrix. This also makes clear that both orientations of a multiple Cartan edge are allowed.

For such a matrix,
\[
 C(H):=H+(2-h)I_n
\]
is a generalized Cartan matrix: its diagonal entries are $2$, its off-diagonal entries are nonpositive integers, and~$c_{ij}=0$ if and only if~$c_{ji}=0$.

\section[Relations between the positivity properties of H and diagonal shifts]{Relations between the positivity properties of $H$ and diagonal shifts $H+\lambda I_n$}\label{s2}
\begin {Lemma} \label{3.1}
For every $H\in M_n(\mathbb R)$ and every scalar $\lambda$,
\[
 \det(H+\lambda I_n)=\sum_{J\subseteq[n]}\lambda^{\,n-|J|}\det H[J],
 \qquad \det H[\varnothing]=1.
\]
\end {Lemma}
\begin{proof}
In the Leibniz expansion of $\det(H+\lambda I_n)$, choose the diagonal term $\lambda$ at precisely the indices outside a subset $J$. The remaining factor is exactly the principal minor $\det H[J]$. Grouping the terms according to $J$ gives the stated identity. \end{proof}

\begin {Proposition} \label{3.2}
If~$H$ is virtual positive semidefinite, then every real eigenvalue of $H$ is nonnegative.
\end {Proposition}
\begin{proof}
Suppose that~$\rho<0$ is a real eigenvalue and put~$\lambda\mathrel{:=}-\rho>0$. Then, $\det(H+\lambda I_n)=0$. On the other hand, Lemma~\ref{3.1} and virtual positive semidefiniteness give
\[
 \det(H+\lambda I_n)
 =\sum_{J\subseteq[n]}\lambda^{n-|J|}\det H[J]
 \ge \lambda^n>0,
\]
a contradiction. \end{proof}

\begin {Lemma} \label{3.111}
\textup{(See \cite[Chapter~4]{Ka90}.)} An indecomposable generalized Cartan matrix is of affine type if and only if it is singular and all its proper principal minors are positive (equivalently, all principal minors are nonnegative and the full determinant is zero).
\end {Lemma}

\begin {Proposition} \label{3.3}
Let~$C$ be a generalized Cartan matrix.

{\rm (i)} $C$ is virtual positive semidefinite if and only if each indecomposable component of $C$ is of finite or affine type.

{\rm (ii)} If~$H\in M_n(\mathbb R)$ is virtual positive semidefinite and~$\lambda>0$, then $H+\lambda I_n$ is virtual positive definite. In particular, if~$H$ is an $h$-Shuhan matrix with~$h<2$ and is virtual positive semidefinite, then $H+(2-h)I_n$ is a generalized Cartan matrix of finite type.

{\rm (iii)} If~$h>2$ and~$H+(2-h)I_n$ is of finite or affine type, then $H$ is virtual positive definite.
\end {Proposition}
\begin{proof}
{\rm (i)} For an indecomposable generalized Cartan matrix, the finite and affine alternatives are characterized by positivity and nonnegativity of the principal minors; see \cite[Chapter~4]{Ka90}. A decomposable matrix is block diagonal after a simultaneous permutation, so the assertion follows componentwise.

{\rm (ii)} Let~$\varnothing\ne J\subseteq[n]$. Since every principal submatrix of a virtual positive semidefinite matrix is again virtual positive semidefinite, Lemma~\ref{3.1}, applied to $H[J]$, gives
\[
 \det\bigl(H[J]+\lambda I_{|J|}\bigr)
 =\sum_{K\subseteq J}\lambda^{|J|-|K|}\det H[K]
 \ge \lambda^{|J|}>0.
\]
Thus, every principal minor of $H+\lambda I_n$ is positive, which proves virtual positive definiteness. If now $H$ is an $h$-Shuhan matrix with~$h<2$, then $C=H+(2-h)I_n$ is a generalized Cartan matrix and is virtual positive definite. Part (i) excludes affine components (their determinants are zero), hence~$C$ is of finite type.

{\rm (iii)} Put~$C\mathrel{:=}H+(2-h)I_n$. By part (i), a finite- or affine-type generalized Cartan matrix is virtual positive semidefinite. Since~$H=C+(h-2)I_n$ and~$h-2>0$, part (ii) shows that~$H$ is virtual positive definite. \end{proof}

\paragraph{Vinberg's trichotomy and spectral criteria.}
An indecomposable generalized Cartan matrix $C$ is said to be of positive, zero, or negative type if there is a vector $u>0$ with $Cu>0$, $Cu=0$, or $Cu<0$, respectively, where the matrix inequalities are understood componentwise \cite{Vin71}.

For an indecomposable $h$-Shuhan matrix of order $n$, put
\[
 B:=hI_n-H=2I_n-C(H),\qquad B_{\rm s}:=\frac{B+B^{\mathsf T}}2.
\]
The entries of both $B$ and $B_{\rm s}$ are nonnegative, and $B$ is irreducible for $n>1$. Writing $\rho$ for spectral radius, we have
\begin{equation}\tag{*}\label{eq:two-spectral-criteria}
\begin{aligned}
 H\text{ is virtual positive semidefinite}&\iff h\ge\rho(B),\\
 H\text{ is generalized positive semidefinite}&\iff h\ge\rho(B_{\rm s}).
\end{aligned}
\end{equation}
To prove the first criterion of~\eqref{eq:two-spectral-criteria}, note that $h<\rho(B)$ gives a negative real eigenvalue of $H$, contradicting Proposition~\ref{3.2}. Conversely, if $h\ge\rho(B)$, adjoin $n-|J|$ zero rows and columns to $B[J]$, supported on $J$, to obtain $0\le\widetilde B_J\le B$. Monotonicity of spectral radius \cite[Chapter~8]{HJ13} gives
\[
 \rho(B[J])\le\rho(B)\le h.
\]
\par\indent For $t>h$, the real factors of $\det(tI_{|J|}-B[J])$ are positive and the nonreal factors occur in complex-conjugate pairs, so this determinant is positive. By continuity, $\det H[J]=\det(hI_{|J|}-B[J])\ge0$. Thus, all principal minors of $H$ are nonnegative.

To prove the second criterion of~\eqref{eq:two-spectral-criteria}, consider $\widehat H=hI_n-B_{\rm s}$. Recall that the spectral theorem for real symmetric matrices states that there are an orthogonal matrix $U$ and real numbers $\lambda_1,\ldots,\lambda_n$ such that
\[
 U^{\mathsf T}B_{\rm s}U=\operatorname{diag}(\lambda_1,\ldots,\lambda_n).
\]
Consequently, $U^{\mathsf T}\widehat H U=\operatorname{diag}(h-\lambda_1,\ldots,h-\lambda_n)$. Since $B_{\rm s}$ is nonnegative and symmetric, $\lambda_{\max}(B_{\rm s})=\rho(B_{\rm s})$ \cite[Chapters~4 and~8]{HJ13}. Thus, $\widehat H$ is positive semidefinite if and only if $h\ge\rho(B_{\rm s})$. By Lemma~\ref{1.2}, this proves the second criterion of~\eqref{eq:two-spectral-criteria}. The corresponding definiteness criteria are $h>\rho(B)$ and $h>\rho(B_{\rm s})$, respectively. When $h=\rho(B)$, the matrix $H=hI_n-B$ is singular; when $h=\rho(B_{\rm s})$, the symmetric part $\widehat H=hI_n-B_{\rm s}$ is singular.

Choose a positive Perron eigenvector $v$ of $B$ (take $v=1$ when $n=1$). Then,
\begin{equation*}\tag{P}\label{eq:perron-shift}
 C(H)v=(2-\rho(B))v.
\end{equation*}
Thus, the three types correspond to $\rho(B)<2$, $\rho(B)=2$, $\rho(B)>2$, respectively. If $h<2$ and $H$ is virtual positive semidefinite, then $\rho(B)\le h<2$ by the first criterion of~\eqref{eq:two-spectral-criteria}. Since $v>0$, \eqref{eq:perron-shift} gives $C(H)v=(2-\rho(B))v>0$. Thus, $C(H)$ is of positive type and hence of finite type; see \cite{Vin71,Ka90}.

The condition $x^{\mathsf T}Hx\ge0$ for every real vector~$x$ can impose a larger threshold. For
\[
 H_h=\begin{pmatrix}h&-2\\-1&h\end{pmatrix},
\]
the virtual and generalized thresholds are $\sqrt2$ and $\tfrac32$. Similarity with respect to the diagonal matrix $D=\operatorname{diag}(\sqrt2,1)$ makes $D^{-1}H_hD$ symmetric; at $h=\sqrt2$ it is positive semidefinite, whereas $H_h$ is not generalized positive semidefinite.

Scalar rescaling to diagonal $2$ need not satisfy Vinberg's product condition. For $H=\left(\begin{smallmatrix}\tfrac32&-1\\-1&\tfrac32\end{smallmatrix}\right)$, the rescaled off-diagonal entries are $-\tfrac43$, whose product $\tfrac{16}{9}$ lies between the consecutive allowed values $1$ and $2$.

\section[Semidefinite h-Shuhan matrices for 0 <= h < 2]{Semidefinite $h$-Shuhan matrices for~$0\le h<2$}\label{s3}
The main classification concerns $0\le h<2$. We also record determinant identities and threshold statements for other values of $h$, with their ranges specified below. Throughout, $H=(h_{ij})$ denotes an indecomposable $h$-Shuhan matrix and~$h\ge0$ unless otherwise stated. For a Cartan matrix~$S$ of finite or affine type, define
\[
  S\psup{h}:=S+(h-2)I,
  \qquad s\psup{h}:=\det S\psup{h},
  \qquad \widehat{s}\psup{h}:=\det\widehat{S\psup{h}}.
\]
Thus, $S\psup{2}=S$, and only the diagonal entries are changed.

\paragraph{Finite-type families.}
We make the matrix conventions explicit. Given a finite Dynkin graph, put~$h$ on every diagonal position, put~$0$ between nonadjacent vertices, and put~$-1$ in both directions on every single edge. On a double (respectively triple) edge, the two directed entries are $-2,-1$ (respectively $-3,-1$). The orientation is part of the matrix data; Definition~\ref{1.7} allows either orientation, and transpose reverses all such arrows. With this convention:
\begin{itemize}
\item $A_n\psup{h}$ is the path on $n$ vertices with all edges single;
\item $B_n\psup{h}$ has $h$ on the diagonal, $(B_n\psup{h})_{i,i+1}=(B_n\psup{h})_{i+1,i}=-1$ for~$1\le i\le n-2$, and at the final edge, $(B_n\psup{h})_{n-1,n}=-2$, $(B_n\psup{h})_{n,n-1}=-1$. We set $C_n\psup{h}:=(B_n\psup{h})^{\mathsf T}$. For~$n\ge3$, these two matrices are generally not permutation equivalent, although all three positivity notions considered here are invariant under transpose;
\item $D_n\psup{h}$ is the simply-laced tree with one trivalent vertex and arms of lengths $n-3,1,1$;
\item $E_6\psup{h},E_7\psup{h},E_8\psup{h}$ are the simply-laced trivalent trees. Their arm lengths are $(2,2,1)$, $(3,2,1)$, and~$(4,2,1)$, respectively;
\item $F_4\psup{h}$ is the four-vertex path $1-2\Rightarrow3-4$, with~$(F_4\psup{h})_{2,3}=-2$, $(F_4\psup{h})_{3,2}=-1$ and all single-edge entries equal to $-1$ in both directions;
\item $G_2\psup{h}=\left(\begin{smallmatrix}h&-3\\-1&h\end{smallmatrix}\right)$.
\end{itemize}
Equivalently, these are obtained from the standard finite-type Cartan matrices in \cite[Chapter~4]{Ka90} by replacing every diagonal $2$ by~$h$. This explicit convention will also be used when a principal submatrix decomposes into smaller Dynkin components.

Set $a_{-1}\psup{h}:=0$, $a_0\psup{h}:=1$, $b_1\psup{h}:=h$, and~$d_3\psup{h}:=a_3\psup{h}$. We also use $\widehat b_0\psup{h}:=1$ and~$\widehat b_1\psup{h}:=h$. The auxiliary value $a_{-1}\psup{h}$ extends the $a$-recurrence to its first step and will also be used when~$i=1$ in Theorem~\ref{5.6}. Expanding the corresponding determinants at an end vertex gives
\begin{equation}\label{e1}
\begin{aligned}
 a_1\psup{h}&=h, & a_n\psup{h}&=h a_{n-1}\psup{h}-a_{n-2}\psup{h} &&(n\ge2),\\
 b_2\psup{h}&=h^2-2, & b_3\psup{h}&=h(h^2-3), & b_n\psup{h}&=h b_{n-1}\psup{h}-b_{n-2}\psup{h} &&(n\ge3),\\
 d_4\psup{h}&=h^2(h^2-3), &&& d_n\psup{h}&=h d_{n-1}\psup{h}-d_{n-2}\psup{h} &&(n\ge5).
\end{aligned}
\end{equation}

The same recurrences also give
\begin{equation}\label{e3}
\begin{aligned}
 b_n\psup{h}&=h a_{n-1}\psup{h}-2a_{n-2}\psup{h}=a_n\psup{h}-a_{n-2}\psup{h} &&(n\ge3),\\
 d_n\psup{h}&=h\bigl(a_{n-1}\psup{h}-a_{n-3}\psup{h}\bigr)=h b_{n-1}\psup{h} &&(n\ge4).
\end{aligned}
\end{equation}
These identities will be used repeatedly below.

Let~$x_n$ satisfy $x_n=hx_{n-1}-x_{n-2}$ whenever the three terms are defined. For~$h\ge0$, $h\ne2$, define $\delta_h$ by
\[
 \delta_h:=
 \begin{cases}
  \sqrt{h^2-4}>0,& h>2,\\
  \mathrm i\sqrt{4-h^2},& 0\le h<2,
 \end{cases}
 \qquad
 \alpha_h:=\frac{h+\delta_h}{2},\quad
 \beta_h:=\frac{h-\delta_h}{2}.
\]
Thus, $\alpha_h$ and~$\beta_h$ are the two distinct roots of $t^2-ht+1=0$; for~$h>2$ they satisfy $\alpha_h>1>\beta_h>0$, while for~$0\le h<2$ they are complex conjugates. Then, \begin {equation}\label{e2}
 x_n=k_x\alpha_h^n+l_x\beta_h^n,
\end {equation}
where~$k_x$ and~$l_x$ are uniquely determined by two consecutive initial values (for example, by~$x_0$ and~$x_1$ when those are specified). We have $\alpha_h+\beta_h=h$, $\alpha_h\beta_h=1$, and~$\alpha_h-\beta_h=\delta_h$. In particular,
$\alpha_h^n=a_{n-1}\psup{h}\alpha_h-a_{n-2}\psup{h}$ and
$\beta_h^n=a_{n-1}\psup{h}\beta_h-a_{n-2}\psup{h}$.
Set
\[
 \Phi:=\frac{1+\sqrt5}{2},\qquad
 \varphi:=\frac{\sqrt5-1}{2}.
\]
Then, $\alpha_{\sqrt5}=\Phi$, $\beta_{\sqrt5}=\varphi$, $\alpha_3=\Phi^2$, and~$\beta_3=\varphi^2$.

For every sequence $x=(x_i)_i$ satisfying the recurrence, define
\[
 \mathcal F_x(i):=h x_i x_{i-1}-x_i^2-x_{i-1}^2.
\]
A substitution of $x_{i+1}=hx_i-x_{i-1}$ shows that~$\mathcal F_x(i+1)=\mathcal F_x(i)$; hence, this quantity is independent of $i$. When the underlying recurrence sequence is clear, we denote this constant by $f(x)$, omitting the index.

\Needspace{7\baselineskip}
\begin{Proposition}\label{4.01}
Assume $h\ge0$ and~$h\ne2$. Then:
\begin{parts}
\item
\[
 a_n\psup{h}=\frac{\alpha_h^{n+1}-\beta_h^{n+1}}{\alpha_h-\beta_h},
 \qquad f(a\psup{h})=-1 \qquad(n\ge1).
\]
\item
\[
 b_n\psup{h}=\alpha_h^n+\beta_h^n,
 \qquad f(b\psup{h})=(\alpha_h-\beta_h)^2 \qquad(n\ge2).
\]
\item
\[
 d_n\psup{h}=h\bigl(\alpha_h^{n-1}+\beta_h^{n-1}\bigr),
 \qquad f(d\psup{h})=(\alpha_h^2-\beta_h^2)^2 \qquad(n\ge4).
\]
\end{parts}
\end{Proposition}
\begin{proof}
For~$a_n\psup{h}$, the initial values $a_0\psup{h}=1$ and~$a_1\psup{h}=h$ in \eqref{e2} give
\[
 k_a=\frac{\alpha_h}{\alpha_h-\beta_h},
 \qquad
 l_a=-\frac{\beta_h}{\alpha_h-\beta_h}.
\]
This proves the formula in part~(i), and the recurrence invariant equals
$k_a l_a(\alpha_h-\beta_h)^2=-1$.

For~$b_n\psup{h}$, the initial values $b_2\psup{h}=h^2-2$ and~$b_3\psup{h}=h(h^2-3)$ give $k_b=l_b=1$. Hence, $b_n\psup{h}=\alpha_h^n+\beta_h^n$ and~$f(b\psup{h})=(\alpha_h-\beta_h)^2$.

Finally, $d_3\psup{h}=h(h^2-2)$ and~$d_4\psup{h}=h^2(h^2-3)$ give
$k_d=h\beta_h$ and~$l_d=h\alpha_h$. Therefore, \[
 d_n\psup{h}=h\bigl(\alpha_h^{n-1}+\beta_h^{n-1}\bigr),
 \qquad
 f(d\psup{h})=h^2(\alpha_h-\beta_h)^2=(\alpha_h^2-\beta_h^2)^2,
\]
where the last equality uses $h=\alpha_h+\beta_h$.
\end{proof}

\Needspace{12\baselineskip}
\begin{Proposition}\label{4.001}
Assume $h\ge0$, $h\ne2$, and~$n\ge1$. Then, we have
\begin{parts}
\item
\begin{align*}
 a_n\psup{h}
 &=\frac1{2^n}\sum_{0\le k<\tfrac{n+1}{2}}
   \binom{n+1}{2k+1}h^{n-2k}(h^2-4)^k,\\
 b_{n+1}\psup{h}
 &=\frac1{2^n}\sum_{0\le k\le\left\lfloor\tfrac{n+1}{2}\right\rfloor}
   \binom{n+1}{2k}h^{n+1-2k}(h^2-4)^k,\\
 d_{n+3}\psup{h}=h b_{n+2}\psup{h}
 &=\frac1{2^{n+1}}\sum_{0\le k\le\left\lfloor\tfrac{n+2}{2}\right\rfloor}
   \binom{n+2}{2k}h^{n+3-2k}(h^2-4)^k.
\end{align*}
\item For~$1\le m<n$, we have
\[
 a_n\psup{h}=a_m\psup{h} a_{n-m}\psup{h}-a_{m-1}\psup{h} a_{n-m-1}\psup{h}.
\]
\item We have
\[
 (b_{n+1}\psup{h})^2-(h^2-4)(a_n\psup{h})^2=4.
\]
\end{parts}
\end{Proposition}
\begin{proof}

Part~(i) follows by applying the binomial theorem to the formulas for $a_n\psup{h}$ and $b_n\psup{h}$ given in parts~(i) and~(ii) of Proposition~\ref{4.01}, together with the identity $d_{n+3}\psup{h}=h b_{n+2}\psup{h}$.

For part~(ii), fix $m$. Both sides satisfy the recurrence $u_n=hu_{n-1}-u_{n-2}$ and agree for $n=m+1$ and $n=m+2$. Hence, the identity holds for every $n>m$.

For part~(iii), substitute the formulas for $a_n\psup{h}$ and $b_n\psup{h}$ from parts~(i) and~(ii) of Proposition~\ref{4.01} into the left-hand side to obtain
\[
 (b_{n+1}\psup{h})^2-(h^2-4)(a_n\psup{h})^2
 =4(\alpha_h\beta_h)^{n+1}=4.
\]
The final equality follows from $\alpha_h\beta_h=1$.

\end{proof}

For later use, note that if~$h>2$, then $\alpha_h>1>\beta_h>0$. Proposition~\ref{4.01} therefore implies positivity of $a_n\psup{h}$, $b_n\psup{h}$, and~$d_n\psup{h}$ in their natural ranges, and
\[
 \frac1{a_n\psup{h}}\longrightarrow0,
 \qquad
 \frac1{b_n\psup{h}}\longrightarrow0,
 \qquad
 \frac1{d_n\psup{h}}\longrightarrow0
 \qquad (n\to\infty).
\]

\begin {Lemma} \label{4.2}
Assume $0\le h<2$, and let~$\theta\in(0,\tfrac{\pi}{2}]$ be the unique number satisfying $\cos\theta=\tfrac{h}{2}$. Then, {\rm (i)}
\[
 a_n\psup{h}=\frac{\sin((n+1)\theta)}{\sin\theta},\qquad
 b_n\psup{h}=2\cos(n\theta),\qquad
 d_n\psup{h}=4\cos\theta\cos((n-1)\theta).
\]

{\rm (ii)} For~$n\ge1$, $a_k\psup{h}\ge0$ for~$1\le k\le n$ if and only if~$h\ge2\cos\frac{\pi}{n+1}$.

{\rm (iii)} For~$n\ge2$, $b_k\psup{h}\ge0$ for~$2\le k\le n$ if and only if~$h\ge2\cos\frac{\pi}{2n}$.

{\rm (iv)} If~$0<h<2$ and~$n\ge4$, then $d_k\psup{h}\ge0$ for~$4\le k\le n$ if and only if~$h\ge2\cos\frac{\pi}{2(n-1)}$.
\end {Lemma}
\begin{proof}
Since~$h=2\cos\theta$, the characteristic roots are $e^{\pm\mathrm{i}\theta}$. Proposition~\ref{4.01}(i), or the elementary solution of the recurrence, gives
$a_n\psup{h}=\tfrac{\sin((n+1)\theta)}{\sin\theta}$. Equation~\eqref{e3} then yields
\[
 b_n\psup{h}=\frac{\sin((n+1)\theta)-\sin((n-1)\theta)}{\sin\theta}=2\cos(n\theta),
\]
and~$d_n\psup{h}=hb_{n-1}\psup{h}=4\cos\theta\cos((n-1)\theta)$.

For (ii), if~$\theta\le\tfrac{\pi}{n+1}$ then $0<(k+1)\theta\le\pi$ for~$1\le k\le n$, so every $a_k\psup{h}\ge0$. Conversely, if~$\theta>\tfrac{\pi}{n+1}$, choose the least integer $r\in\{2,\ldots,n+1\}$ with~$r\theta>\pi$. Then, $(r-1)\theta\le\pi$ and, because $\theta\le\tfrac{\pi}{2}$, one has $\pi<r\theta\le\tfrac{3\pi}{2}$, so~$a_{r-1}\psup{h}<0$. Thus, $\theta\le\tfrac{\pi}{n+1}$, equivalently $h\ge2\cos\!\left(\tfrac{\pi}{n+1}\right)$.

For (iii), the same first-crossing argument applied to $b_k\psup{h}=2\cos(k\theta)$ shows that all $b_k\psup{h}$ for~$2\le k\le n$ are nonnegative exactly when~$n\theta\le\tfrac{\pi}{2}$, i.e. $h\ge2\cos\!\left(\tfrac{\pi}{2n}\right)$. For (iv), $d_4\psup{h}=h^2(h^2-3)\ge0$ and~$h>0$ first give $h\ge\sqrt3$, hence~$b_2\psup{h}>0$. The remaining conditions $d_k\psup{h}\ge0$ are equivalent to $b_{k-1}\psup{h}\ge0$. Part (iii), now applied to the full range $2\le k-1\le n-1$, gives the asserted threshold. \end{proof}

\begin {Lemma}\label{finite-principal-structure}
Every principal submatrix of a finite-type family is, after a simultaneous permutation, a direct sum of connected blocks of the following kinds:
\begin{enumerate}
\item a principal submatrix of $A_n\psup{h}$ is a direct sum of $A_k\psup{h}$ blocks;
\item a principal submatrix of $B_n\psup{h}$ (respectively $C_n\psup{h}$) is a direct sum of $A_k\psup{h}$ blocks and at most one $B_\ell\psup{h}$ (respectively $C_\ell\psup{h}$) block;
\item a principal submatrix of $D_n\psup{h}$ is a direct sum of $A_k\psup{h}$ and (possibly) one smaller $D_\ell\psup{h}$ block;
\item a principal submatrix of $E_m\psup{h}$ ($m=6,7,8$) is a direct sum of blocks of type $A$, $D$, and smaller $E$ whenever the latter occurs as a connected induced subdiagram;
\item a proper principal submatrix of $F_4\psup{h}$ is a direct sum of $A$-type blocks and at most one multiple-edge block, of type $B_2\psup{h}$ or (in rank three) $B_3\psup{h}$ or $C_3\psup{h}$; for~$G_2\psup{h}$ the only proper connected principal blocks are $A_1\psup{h}$.
\end{enumerate}
\end {Lemma}
\begin{proof}
The Dynkin graphs listed above are trees. Deleting vertices produces an induced forest, and after simultaneously permuting rows and columns according to its connected components the principal submatrix becomes block diagonal. A connected component not containing a multiple edge is a simply-laced path or subtree of the indicated $A,D,E$ type. Since the $B_n,C_n$ and~$F_4$ diagrams have only one multiple edge, at most one component contains it, producing the stated oriented multiple-edge block. The $G_2$ assertion is immediate. \end{proof}

\begin {Proposition} \label{4.3}
{\rm (i)} $A_{n}\psup{h}$ is positive semidefinite if and only if~$h\geq 2\cos\frac{\pi}{n+1}$, $n\geq1$.

{\rm (ii)} $B_{n}\psup{h}$ and~$C_n\psup{h}=(B_n\psup{h})^{\mathsf T}$ are virtual positive semidefinite if and only if~$h\geq 2\cos\frac{\pi}{2n}$, $n\geq2$.

{\rm (iii)} $D_{n}\psup{h}$ is positive semidefinite if and only if~$h\geq 2\cos\frac{\pi}{2(n-1)}$, $n\geq4$.

{\rm (iv)} $E_{6}\psup{h}$ is positive semidefinite if and only if~$h\geq2\cos\frac{\pi}{12}$.

{\rm (v)} $E_{7}\psup{h}$ is positive semidefinite if and only if~$h\geq2\cos\frac{\pi}{18}$.

{\rm (vi)} $E_{8}\psup{h}$ is positive semidefinite if and only if~$h\geq2\cos\frac{\pi}{30}$.

{\rm (vii)} $F_{4}\psup{h}$ is virtual positive semidefinite if and only if~$h\geq2\cos\frac{\pi}{12}$.

{\rm (viii)} $G_{2}\psup{h}$ is virtual positive semidefinite if and only if~$h\geq2\cos\frac{\pi}{6}$.
\end {Proposition}
\begin{proof}

For $h\ge2$, the matrices obtained from the finite-type Cartan matrices by replacing each diagonal entry $2$ by $h$ have strictly positive principal minors. At $h=2$ this follows from the principal-minor characterization of finite-type Cartan matrices (see~\cite[Chapter~4]{Ka90}); for $h>2$ it follows from Proposition~\ref{3.3}(ii). Consequently, the matrices in the symmetric $A,D,E$ families are positive definite. At $h=0$, each matrix in the finite-type families under consideration, except $A_1\psup{0}$, contains a principal $2\times2$ submatrix corresponding to two adjacent vertices of its Dynkin diagram. This submatrix has the form $\left(\begin{smallmatrix}0&-a\\-b&0\end{smallmatrix}\right)$ with $a,b>0$, so its determinant is $-ab<0$. Thus, among the displayed finite-type families, only the $1\times1$ zero matrix $A_1\psup{0}$ is positive semidefinite. It remains to consider $0<h<2$.

By Lemma~\ref{finite-principal-structure}, virtual positive semidefiniteness is equivalent to nonnegativity of the determinants of the connected principal blocks listed there. For the symmetric matrices of types $A$, $D$, and $E$, nonnegativity of all principal minors is also equivalent to positive semidefiniteness.

Parts (i)--(iii) now follow from Lemma~\ref{4.2}: for~$A_n$ the strongest condition is $a_n\psup{h}\ge0$ together with the preceding $a_k\psup{h}$; for~$B_n$ it is the simultaneous nonnegativity of $a_{k-1}\psup{h}$ and~$b_k\psup{h}$, whose sharp maximum threshold is $2\cos\!\left(\tfrac{\pi}{2n}\right)$; the same holds for~$C_n\psup{h}=(B_n\psup{h})^{\mathsf T}$ because principal determinants are unchanged by transpose; and for~$D_n$ it is the simultaneous nonnegativity of the relevant $A$-blocks and~$d_k\psup{h}$, giving $2\cos\!\left(\tfrac{\pi}{2(n-1)}\right)$.

\Needspace{8\baselineskip}
For the exceptional types, the structural lemma reduces the verification to the full determinant plus already controlled proper connected subdiagrams. Direct end-vertex expansions give
\[
\begin{aligned}
 e_6\psup{h}&=(h^2-1)\bigl(h^2-(2+\sqrt3)\bigr)\bigl(h^2-(2-\sqrt3)\bigr),\\
 f_4\psup{h}&=\bigl(h^2-(2+\sqrt3)\bigr)\bigl(h^2-(2-\sqrt3)\bigr),\\
 g_2\psup{h}&=h^2-3.
\end{aligned}
\]
Together with the proper-subdiagram thresholds, these give (iv), (vii), and (viii). For~$E_7$, writing $h=2\cos\theta$ gives
$e_7\psup{h}=2\cos\theta\,(2\cos6\theta-1)$, so on the range already forced by the $E_6$ subdiagram the sharp zero is $\theta=\tfrac{\pi}{18}$. For~$E_8$, the $E_7$ condition restricts $0\le\theta\le\tfrac{\pi}{18}$ and
\[
 e_8\psup{h}=2\{-2\sin(5\theta)\sin(3\theta)+\cos(6\theta)-\tfrac12\}.
\]
On this interval, the expression is strictly decreasing as a function of $\theta$: its derivative is

\[
\begin{aligned}
 \frac{d}{d\theta}e_8\psup{2\cos\theta}
 &=-4\bigl[5\cos(5\theta)\sin(3\theta)
             +3\sin(5\theta)\cos(3\theta)\bigr]\\
 &\hspace{1.5em}-12\sin(6\theta)<0
 \qquad(0<\theta\le\tfrac{\pi}{18}).
\end{aligned}
\]

The determinant vanishes at~$\theta=\tfrac{\pi}{30}$, yielding (vi). Every proper connected principal subdiagram has a threshold no greater than $2\cos(\tfrac{\pi}{30})$, so no further condition is needed. \end{proof}

To study the determinant polynomials $\widehat b_n\psup{h}$ uniformly in $n$,
we place the exceptional edge of weight $\tfrac32$ between the first two
coordinates. More precisely, for $n\ge2$, let $J_n=(j_{ij})$ be the
real symmetric tridiagonal matrix with zero diagonal, whose only
nonzero entries are
\[
 j_{i,i+1}=j_{i+1,i}=
 \begin{cases}
  \tfrac32,&i=1,\\
  1,&2\le i\le n-1.
 \end{cases}
\]
Thus,
\[
 J_2=\begin{pmatrix}0&\tfrac32\\[1pt]\tfrac32&0\end{pmatrix},
 \qquad
 J_3=\begin{pmatrix}0&\tfrac32&0\\[1pt]\tfrac32&0&1\\[1pt]0&1&0\end{pmatrix}.
\]
Let $P_n$ denote the permutation matrix reversing coordinate order:
$P_ne_i=e_{n+1-i}$. In the symmetric part of $B_n\psup{h}$, the last
pair of adjacent coordinates has entries $-\tfrac32$ and all other
adjacent pairs have entries $-1$. Consequently,
\begin{equation}\label{eq:Jn-reversal}
 P_n^{\mathsf T}\widehat{B_n\psup{h}}P_n=hI_n-J_n,
 \qquad \widehat b_n\psup{h}=\det(hI_n-J_n).
\end{equation}

\begin {Lemma} \label{4.5''}
All roots of $\widehat b_n\psup{h}=\det\widehat{B_n\psup{h}}$, regarded as a polynomial in $h$, are real and simple.
\end {Lemma}
\begin{proof}

By \eqref{eq:Jn-reversal}, the roots are the eigenvalues of the real
symmetric matrix $J_n$, and hence are real. To prove simplicity, let
$J_nv=\lambda v$ with $v=(v_1,\ldots,v_n)^{\mathsf T}\ne0$.
The first row gives $\tfrac32v_2=\lambda v_1$; for
$2\le k\le n-1$, the $k$th row determines $v_{k+1}$ uniquely from
$v_{k-1}$ and $v_k$, because $j_{k,k+1}=1$. Thus $v_1=0$ would imply
$v=0$, and every eigenvector is determined by its first coordinate.
Each eigenspace is therefore one-dimensional. Since a real symmetric
matrix is diagonalizable, all its eigenvalues are simple.

\end{proof}

Recall the following standard terminology and results. For an irreducible nonnegative matrix $M$, a \emph{Perron vector} is a vector $u>0$ such that $Mu=\rho(M)u$, where $\rho(M)$ is the spectral radius. The Perron--Frobenius theorem asserts that such a vector exists and that $\rho(M)$ is an algebraically simple eigenvalue. For a real symmetric matrix $A$ and $x\ne0$, the \emph{Rayleigh quotient} is $R_A(x):=x^{\mathsf T}Ax/(x^{\mathsf T}x)$. The Rayleigh--Ritz variational principle states that $\lambda_{\max}(A)=\max_{x\ne0}R_A(x)$. See \cite[Chapters~4 and~8]{HJ13}.

Write $\mu_n$ for the largest real zero of $\widehat b_n\psup{h}$.
By \eqref{eq:Jn-reversal}, $\mu_n=\lambda_{\max}(J_n)$.
Moreover, $J_n$ is entrywise nonnegative and irreducible, so
$\mu_n=\rho(J_n)>0$ by the Perron--Frobenius theorem.

\begin {Lemma} \label{4.6}
For every $n\ge3$, the largest roots satisfy $\mu_n>\mu_{n-1}$.
\end {Lemma}
\begin{proof}

For $n\ge3$, the definition of $J_n$ shows that its leading principal
submatrix of order $n-1$ is exactly $J_{n-1}$. The only nonzero entry
connecting this block to the last coordinate is
$(J_n)_{n-1,n}=1$. Thus,
\[
 J_n=\begin{pmatrix}
        J_{n-1}&e_{n-1}\\
        e_{n-1}^{\mathsf T}&0
       \end{pmatrix},
\]
where $e_{n-1}=(0,\ldots,0,1)^{\mathsf T}\in\mathbb R^{n-1}$.

\noindent By the Perron--Frobenius theorem, there is a vector
$u=(u_1,\ldots,u_{n-1})^{\mathsf T}$ with every $u_i>0$ such that
$J_{n-1}u=\mu_{n-1}u$.

For $t>0$, let $z_t=(u^{\mathsf T},t)^{\mathsf T}$.
Multiplying out the displayed block matrix, we obtain
\[
\begin{aligned}
 z_t^{\mathsf T}J_nz_t
 &=u^{\mathsf T}J_{n-1}u
   +t\,u^{\mathsf T}e_{n-1}
   +t\,e_{n-1}^{\mathsf T}u\\
 &=\mu_{n-1}\lVert u\rVert^2+2tu_{n-1},\\
 z_t^{\mathsf T}z_t&=\lVert u\rVert^2+t^2.
\end{aligned}
\]
It follows that the Rayleigh quotient
$R_{J_n}(z_t)=z_t^{\mathsf T}J_nz_t/(z_t^{\mathsf T}z_t)$ satisfies
\[
\begin{aligned}
 R_{J_n}(z_t)-\mu_{n-1}
 &=\frac{\mu_{n-1}\lVert u\rVert^2+2tu_{n-1}}
         {\lVert u\rVert^2+t^2}-\mu_{n-1}\\
 &=\frac{t(2u_{n-1}-\mu_{n-1}t)}{\lVert u\rVert^2+t^2}>0
 \qquad\left(0<t<\frac{2u_{n-1}}{\mu_{n-1}}\right).
\end{aligned}
\]
Finally, the Rayleigh--Ritz characterization of the largest eigenvalue
of a real symmetric matrix gives
\[
 \mu_n=\lambda_{\max}(J_n)
 =\max_{z\ne0}R_{J_n}(z)
 \ge R_{J_n}(z_t)>\mu_{n-1}.
\]

\end{proof}

Let~$\epsilon$ be the largest real root of
\[
 g(h):=\widehat b_3\psup{h}-\widehat b_2\psup{h}
      =h^3-h^2-\frac{13}{4}h+\frac94.
\]
Numerically, $\epsilon\approx2.04998$.

\begin {Lemma} \label{4.7} The bound $\tfrac32\le\mu_n<\epsilon$ holds for every $n\ge2$.
\end {Lemma}
\begin{proof}

The lower bound follows from $\mu_2=\tfrac32$ and Lemma~\ref{4.6}. Since~$g(2)=-\tfrac14$, $g'(2)=\tfrac{19}{4}$, and~$g''(h)=6h-2>0$ for~$h\ge2$, the largest root $\epsilon$ is greater than $2$ and~$g(h)\ge0$ for~$h\ge\epsilon$. The initial thresholds satisfy
\[
 \mu_2=\tfrac32<2<\epsilon,\qquad
 \mu_3=\tfrac{\sqrt{13}}2<2<\epsilon.
\]
Fix $h\ge\epsilon$ and put~$q_k:=\widehat b_k\psup{h}$. Then, $q_2,q_3>0$ and~$q_3-q_2=g(h)\ge0$. For~$k\ge4$, the determinant recurrence gives
\[
 q_k-q_{k-1}=(h-2)q_{k-1}+(q_{k-1}-q_{k-2}).
\]
Induction yields $q_k>q_{k-1}>0$ for every $k\ge4$. Thus, none of the polynomials $\widehat b_k\psup{h}$ has a zero in $[\epsilon,\infty)$, and~$\mu_k<\epsilon$ for all $k\ge2$.

\end{proof}

Equivalently, $\mu_n=\rho(J_n)$ is the spectral threshold for generalized positive definiteness of both $B_n\psup{h}$ and~$C_n\psup{h}=(B_n\psup{h})^{\mathsf T}$. Lemmas~\ref{4.6} and \ref{4.7} show that these thresholds increase strictly with~$n$ and remain below $\epsilon\approx2.04998$.

\begin {Corollary} \label{4.18} If~$h\geq\epsilon$, then both $B_{n}\psup{h}$ and~$C_n\psup{h}$ are generalized positive definite for every $n\geq2$.
\end {Corollary}

For the $B_n,C_n$ families, the restriction $h<2$ requires $\mu_n<2$. We show that $\mu_9=2$. At $h=2$, the determinant recurrence with initial values $2,\tfrac74$ gives
\[
 \widehat b_k\psup{2}=\frac{9-k}{4}\qquad(k\ge1).
\]
Set
\[
 Q:=2I_9-J_9=\begin{pmatrix}Q_8&z\\z^{\mathsf T}&2\end{pmatrix}.
\]
The leading principal minors of $Q_8$ are
$\widehat b_k\psup{2}=(9-k)/4>0$ for $1\le k\le8$.
Hence $Q_8$ is positive definite by Sylvester's criterion.
\par\Needspace{12\baselineskip}
Recall the Schur complement criterion. Let
$M=\left(\begin{smallmatrix}A&z\\z^{\mathsf T}&d\end{smallmatrix}\right)$
be a real symmetric block matrix with $A$ positive definite, and put
$s=d-z^{\mathsf T}A^{-1}z$, the Schur complement of $A$.
Block elimination gives the congruence

\[
\begin{pmatrix}I&0\\-z^{\mathsf T}A^{-1}&1\end{pmatrix}
M
\begin{pmatrix}I&-A^{-1}z\\0&1\end{pmatrix}
=\begin{pmatrix}A&0\\0&s\end{pmatrix}.
\]

Since the two outer matrices are invertible transposes of one another,
$M$ is positive semidefinite if and only if $s\ge0$; moreover,
$\det M=(\det A)s$ (see~\cite[Chapter~7]{HJ13}).
For $Q$, the Schur complement of $Q_8$ is
$s=2-z^{\mathsf T}Q_8^{-1}z$.
As $\det Q=\widehat b_9\psup{2}=0$ and $\det Q_8>0$,
the determinant formula gives $s=0$.
Thus $Q=2I_9-J_9$ is positive semidefinite and singular,
so the largest eigenvalue of $J_9$ is $2$, that is, $\mu_9=2$.
Lemmas~\ref{4.6} and~\ref{4.7} now yield

\[
 \frac32=\mu_2<\cdots<\mu_8<\mu_9=2<\mu_{10}<\cdots<\epsilon.
\]

\begin {Proposition} \label{4.8} {\rm (i)} If~$0\leq h<2$, then $B_{n}\psup{h}$ and~$C_n\psup{h}=(B_n\psup{h})^{\mathsf T}$ are generalized positive semidefinite if and only if~$h\geq\mu_{n}$ for~$2\leq n\leq8$, where
\[
\mu_{n}=\left \{\begin{array} {l}
\frac{3}{2}\ , \ n=2\\[1mm]
\frac{\sqrt{13}}{2}\ ,\ n=3\\[1mm]
\sqrt{\frac{17}{8}+\frac{\sqrt{145}}{8}}\ , \ n=4\\[1mm]
\sqrt{\frac{21}{8}+\frac{\sqrt{89}}{8}}\ , \ n=5\\[1mm]
\sqrt{\frac{25}{12}+\frac{\sqrt{157}}{6}\cos\theta_{1}}\ , \ n=6\\[1mm]
\sqrt{\frac{29}{12}+\frac{11}{12}\cos\theta_{2}+\frac{11\sqrt{3}}{12}
\sin\theta_{2}}\ ,\ n=7\\[1mm]
\sqrt{\frac{33}{16}+\frac{\beta}{16}+
\frac{1}{2}\sqrt{\frac{547}{96}-\frac{\alpha}{6}+\frac{17}{32\beta}}}\ ,
\ n=8.\end {array} \right.
\]

where
\[
\begin{aligned}
 \theta_1&:=\frac13\arctan\frac{54\sqrt{1327}}{19},&
 \theta_2&:=\frac13\arctan\frac{12\sqrt{11919}}{235},\\
 \alpha&:=\sqrt{727}\cos\!\left(\frac13\arctan\frac{3\sqrt{37701987}}{34607}\right),&
 \beta&:=\sqrt{\frac{547}{3}+\frac{32}{3}\alpha}.
\end{aligned}
\]

The corresponding numerical values are

\[
\begin{aligned}
 \mu_2&\approx1.500000,& \mu_3&\approx1.802776,& \mu_4&\approx1.905308,\\
 \mu_5&\approx1.950448,& \mu_6&\approx1.973634,& \mu_7&\approx1.986807,\\
 \mu_8&\approx1.994840.
\end{aligned}
\]

{\rm (ii)} $F_{4}\psup{h}$ is generalized positive semidefinite if and only if~$h\geq2$.

{\rm (iii)} $G_{2}\psup{h}$ is generalized positive semidefinite if and only if~$h\geq2$.
\end {Proposition}
\begin{proof}
By Lemma~\ref{1.2}, $B_n\psup{h}$ is generalized positive semidefinite if and only if
$\widehat{B_n\psup{h}}$ is positive semidefinite. By
\eqref{eq:Jn-reversal}, this is equivalent to the positive semidefiniteness
of $hI_n-J_n$, and hence to
\[
  h\ge \lambda_{\max}(J_n)=\mu_n.
\]
This explains why the \emph{largest} root, rather than an arbitrary root of $\widehat b_n\psup{h}$, is the relevant threshold. Since~$x^{\mathsf T}B_n\psup{h} x=x^{\mathsf T}(B_n\psup{h})^{\mathsf T}x$, exactly the same threshold holds for~$C_n\psup{h}$. Lemmas~\ref{4.6} and \ref{4.7}, together with~$\mu_9=2$, then show that only $2\le n\le8$ can occur in the range $0\le h<2$.

For transparency, the determinant polynomials whose largest roots are used above are
{\renewcommand{\arraystretch}{1.45}
\[
\begin{array}{c|l}
 n & \widehat b_n\psup{h} \\ \hline
2& h^2-\frac94\\
3& h^3-\frac{13}{4}h\\
4& h^4-\frac{17}{4}h^2+\frac94\\
5& h^5-\frac{21}{4}h^3+\frac{11}{2}h\\
6& h^6-\frac{25}{4}h^4+\frac{39}{4}h^2-\frac94\\
7& h^7-\frac{29}{4}h^5+15h^3-\frac{31}{4}h\\
8& h^8-\frac{33}{4}h^6+\frac{85}{4}h^4-\frac{35}{2}h^2+\frac94.
\end{array}
\]}
They are generated recursively from the tridiagonal determinant (equivalently from the recurrence $\widehat b_n\psup{h}=h\widehat b_{n-1}\psup{h}-\widehat b_{n-2}\psup{h}$ once the initial cases are fixed). Thus, the displayed radical and trigonometric expressions for~$\mu_2,\ldots,\mu_8$ can be checked by substitution into these explicit polynomials. Appendix~\ref{app:cas} gives a short exact-arithmetic script reproducing the polynomials and their largest real roots.

For~$F_4\psup{h}$ and~$G_2\psup{h}$, generalized positive semidefiniteness is again equivalent to positive semidefiniteness of the symmetrization. For~$F_4\psup{h}$, the weighted path in the symmetrization has off-diagonal weights $1,\tfrac{3}{2},1$, so
\[
 \det\widehat{F_4\psup{h}}=h^4-\frac{17}{4}h^2+1=\bigl(h^2-4\bigr)\bigl(h^2-\tfrac14\bigr).
\]
Thus, its largest adjacency eigenvalue is $2$. For~$G_2\psup{h}$, the symmetrization is
$\left(\begin{smallmatrix}h&-2\\-2&h\end{smallmatrix}\right)$, whose eigenvalues are $h\pm2$. Hence, both families are generalized positive semidefinite exactly for~$h\ge2$. \end{proof}
\begin {Remark} \label{4.9}
For the classification theorem, only the largest real zero of each polynomial $\widehat b_n\psup{h}$ is needed.  In particular, the value $\mu_8\approx1.994839578989595$ is independently reproduced by the exact-arithmetic computation in Appendix~\ref{app:cas}; no Galois-theoretic assertion is used anywhere in the proof.
\end {Remark}

The reduction from $H$ to the finite-type Cartan matrix $C(H)=H+(2-h)I_n$ in Proposition~\ref{3.3}(ii), together with the preceding spectral thresholds, gives the following classification.

\Needspace{24\baselineskip}
\begin{Theorem}\label{4.4}
Assume $0\le h<2$. Then, the following classifications hold.
\begin{mainparts}
\item $H$ is positive semidefinite if and only if, after a simultaneous permutation of rows and columns, it is one of
\begin{parts}
\item $A_n\psup{h}$ with~$n\ge1$ and~$h\ge2\cos\!\left(\tfrac{\pi}{n+1}\right)$;
\item $D_n\psup{h}$ with~$n\ge4$ and~$h\ge2\cos\!\left(\tfrac{\pi}{2(n-1)}\right)$;
\item $E_6\psup{h},E_7\psup{h},E_8\psup{h}$ with thresholds $2\cos\!\left(\tfrac{\pi}{12}\right)$, $2\cos\!\left(\tfrac{\pi}{18}\right)$, and~$2\cos\!\left(\tfrac{\pi}{30}\right)$, respectively.
\end{parts}

\item $H$ is virtual positive semidefinite if and only if it is positive semidefinite or, after a simultaneous permutation, it is one of
\begin{parts}
\item $B_n\psup{h}$ or $C_n\psup{h}=(B_n\psup{h})^{\mathsf T}$, with~$n\ge2$ and~$h\ge2\cos\!\left(\tfrac{\pi}{2n}\right)$;
\item $F_4\psup{h}$ with~$h\ge2\cos\!\left(\tfrac{\pi}{12}\right)$;
\item $G_2\psup{h}$ with~$h\ge2\cos\!\left(\tfrac{\pi}{6}\right)$.
\end{parts}

\item $H$ is generalized positive semidefinite if and only if it is positive semidefinite or, after a simultaneous permutation,
\[
 H=B_n\psup{h}\quad\text{or}\quad H=C_n\psup{h}=(B_n\psup{h})^{\mathsf T},
\]
where~$2\le n\le8$ and~$h\ge\mu_n$. The thresholds $\mu_n$ are listed in Proposition~\ref{4.8}.
\end{mainparts}
\end{Theorem}
\begin{proof}
We give the details in the three steps on which the classification rests.

\emph{Step 1: reduction to finite type.}
Suppose first that~$H$ is virtual positive semidefinite. This includes the positive and generalized positive semidefinite cases by Proposition~\ref{1.4}. Because $h<2$, Proposition~\ref{3.3}(ii) shows that
$C=H+(2-h)I_n$ is a finite-type Cartan matrix. Since~$H$ is indecomposable, so is $C$.

\emph{Step 2: exhaustion of the finite families.}
The classification of indecomposable finite-type Cartan matrices gives the types
$A_n,B_n,C_n,D_n,E_6,E_7,E_8,F_4,G_2$; see \mbox{\cite[Chapter~4]{Ka90}}. Both orientations of the double edge occur in types $B_n$ and~$C_n$; for~$n\ge3$ they are transposes rather than permutation-equivalent matrices. (For~$F_4$ and~$G_2$, the opposite orientation is permutation equivalent to the displayed representative.) Consequently, $H_\sigma$ must be one of
\[
 A_n\psup{h},\ B_n\psup{h},\ C_n\psup{h},\ D_n\psup{h},\ E_6\psup{h},\ E_7\psup{h},\ E_8\psup{h},\ F_4\psup{h},\ G_2\psup{h}.
\]
There are no further cases. Proposition~\ref{4.3} gives the exact virtual and positive-semidefinite thresholds for these families, while Proposition~\ref{4.8} gives the generalized thresholds for the nonsymmetric families.

\emph{Step 3: separation of the three positivity notions.}
For positive semidefiniteness, the matrix must be symmetric. Among the representatives above this leaves precisely $A,D,E$; applying Proposition~\ref{4.3}(i),(iii)--(vi) gives part (1).

For virtual positive semidefiniteness, all finite types may occur. Proposition~\ref{4.3} gives exactly the additional nonsymmetric families $B_n\psup{h},C_n\psup{h},F_4\psup{h},G_2\psup{h}$ with the thresholds stated in part (2), proving both necessity and sufficiency.

For generalized positive semidefiniteness, the symmetric matrices of types $A$, $D$, and $E$ again give part (1). For~$B_n\psup{h}$ and~$C_n\psup{h}$, Proposition~\ref{4.8}(i) gives $h\ge\mu_n$. Since~$\mu_9=2$ and Lemma~\ref{4.6} gives strict increase, the restriction $h<2$ leaves exactly $2\le n\le8$. Proposition~\ref{4.8}(ii),(iii) shows that~$F_4\psup{h}$ and~$G_2\psup{h}$ would require $h\ge2$, so neither contributes in the present range. This proves part (3). The sufficiency assertions in all three parts are exactly the corresponding propositions already cited. \end{proof}

\paragraph{Affine-type notation.}
Let~$X_N^{(r)}$ denote the standard affine generalized Cartan matrix in the notation of \cite[Chapter~4]{Ka90}, and let~$m$ be its order. Define
\begin{equation}\label{affine-h-notation}
 X_N\affpsup{r}{h}:=X_N^{(r)}+(h-2)I_m.
\end{equation}
Thus, only the diagonal entries $2$ are changed to $h$; all off-diagonal Cartan entries remain unchanged. In the comparisons below, $S$ denotes the relevant unshifted affine matrix. For the nonsymmetric comparisons immediately below, we use the standard affine types

\[
\begin{gathered}
 G_2^{(1)},\ D_4^{(3)},\ F_4^{(1)},\ E_6^{(2)},\ A_2^{(2)},\\
 B_n^{(1)},\ A_{2n-1}^{(2)},\ C_n^{(1)},\ A_{2n}^{(2)},\ D_{n+1}^{(2)},
\end{gathered}
\]

with the standard superscripts and ranks; symmetric affine types used later retain their usual symbols. 
For the ranges used below, the weighted graphs of the symmetrizations $\tfrac12(S+S^{\mathsf T})$ agree, up to relabelling of vertices, in the pairs (or triples)

\[
\begin{gathered}
 (G_2^{(1)},D_4^{(3)}),\qquad (F_4^{(1)},E_6^{(2)}),\\
 (B_n^{(1)},A_{2n-1}^{(2)})\ (n>2),\qquad
 (C_n^{(1)},A_{2n}^{(2)},D_{n+1}^{(2)})\ (n>1).
\end{gathered}
\]

A single edge has symmetrized weight $1$, a double edge weight $\tfrac{3}{2}$, and a triple edge weight $2$; hence, reversing an arrow does not change the weighted graph. The $B_n^{(1)}$ pair gives a fork with one double end-edge and the $C_n^{(1)}$ triple a path with double end-edges; $A_2^{(2)}$ is treated separately. Hence, the paired types have the same generalized-semidefinite threshold, although the original Cartan matrices need not be permutation equivalent.

\begin {Proposition} \label{4.5'}
If~$S$ is an indecomposable affine generalized Cartan matrix, then $S\psup{h}$ is virtual positive semidefinite if and only if~$h\ge2$.
\end {Proposition}
\begin{proof}
If~$h=2$, then $S\psup{2}=S$ is virtual positive semidefinite by Proposition~\ref{3.3}(i). If~$h>2$, then $S\psup{h}=S+(h-2)I$ is virtual positive definite by Proposition~\ref{3.3}(ii). Conversely, if~$h<2$ and~$S\psup{h}$ were virtual positive semidefinite, Proposition~\ref{3.3}(ii) would imply that~$S=S\psup{h}+(2-h)I$ is of finite type, contradicting that~$S$ is affine. \end{proof}

\begin {Proposition} \label{4.5}
{\rm (i)} $G_{2}\affpsup{1}{h}$ (and~$D_{4}\affpsup{3}{h}$) is generalized positive semidefinite if and only if~$h\geq \sqrt{5}$.

{\rm (ii)} $F_{4}\affpsup{1}{h}$ (and~$E_{6}\affpsup{2}{h}$) is generalized positive semidefinite if and only if~$h\geq \frac{\sqrt{17}}{2}$.

{\rm (iii)} $A_{2}\affpsup{2}{h}$ is generalized positive semidefinite if and only if  $h\geq \frac{5}{2}$.

\end {Proposition}
\begin{proof}
For each matrix~$S\psup{h}$ in (i)--(iii), Lemma~\ref{1.2} gives
\[
   S\psup{h}\text{ generalized positive semidefinite}
   \iff \widehat{S\psup{h}}=hI-K\succeq0,
\]
where~$K$ is a real symmetric nonnegative matrix. Hence, the exact threshold is the largest eigenvalue $\rho(K)$, equivalently the largest real zero of $\det\widehat{S\psup{h}}$.
For the three displayed affine types, the determinant computations are

\[
\begin{aligned}
 \widehat g_2\affpsup{1}{h}
 &=\widehat d_4\affpsup{3}{h}=h(h^2-5),\\
 \widehat f_4\affpsup{1}{h}
 &=\widehat e_6\affpsup{2}{h}
   =h(h^2-1)(h^2-\tfrac{17}{4}),\\
 \widehat a_2\affpsup{2}{h}&=h^2-\frac{25}{4}.
\end{aligned}
\]

The largest real zeros are respectively $\sqrt5$, $\tfrac{\sqrt{17}}{2}$, and~$\tfrac{5}{2}$, which proves all three equivalences. \end{proof}

\begin {Lemma} \label{4.10}
For every $n\ge3$ and every $h>\epsilon$,
\[
 \widehat b_n\psup{h}>\widehat b_{n-1}\psup{h}>\cdots>\widehat b_2\psup{h}>0.
\]
\end {Lemma}
\begin{proof}
By Lemma~\ref{4.7}, the largest zero $\mu_k$ of $\widehat b_k\psup{h}$ satisfies $\mu_k<\epsilon$ for every $k$, so~$\widehat b_k\psup{h}>0$ whenever $h>\epsilon$. Moreover,
\[
 \widehat b_k\psup{h}-\widehat b_{k-1}\psup{h}
 =(h-1)\widehat b_{k-1}\psup{h}-\widehat b_{k-2}\psup{h}.
\]
Starting with
$\widehat b_3\psup{h}-\widehat b_2\psup{h}=g(h)>g(\epsilon)=0$
for~$h>\epsilon$, induction gives
$\widehat b_k\psup{h}-\widehat b_{k-1}\psup{h}>0$ for every $k\ge3$. \end{proof}

\begin {Lemma} \label{4.10'}
The largest roots of $\widehat b_n\affpsup{1}{h}$ and~$\widehat c_n\affpsup{1}{h}$, regarded as polynomials in $h$, are simple.
\end {Lemma}
\begin{proof}
For either affine family, the symmetrized matrix has the form
\[
   \widehat{S\psup{h}}=hI-K,
\]
where~$K$ is a real symmetric nonnegative irreducible weighted adjacency matrix of the connected affine Dynkin graph. Hence, the largest zero of $\det(hI-K)$ is the spectral radius $\rho(K)$. By the Perron--Frobenius theorem, $\rho(K)$ is a simple eigenvalue of the irreducible nonnegative matrix~$K$. Therefore, the largest root of the characteristic polynomial is simple. \end{proof}

Define
\[
 \lambda_n:=\max\{h\in\mathbb R:\widehat b_n\affpsup{1}{h}=0\},\qquad
 \eta_n:=\max\{h\in\mathbb R:\widehat c_n\affpsup{1}{h}=0\}.
\]

\begin {Proposition} \label{4.11}
{\rm (i)} For every $n>2$, $B_{n}\affpsup{1}{h}$ (and~$A_{2n-1}\affpsup{2}{h}$) is generalized positive semidefinite whenever $h\geq\frac{\sqrt{17}}{2}$. This uniform bound is sharp over the whole family, since equality is required when~$n=3$.

{\rm (ii)} For every $n>1$, $C_{n}\affpsup{1}{h}$ (and~$A_{2n}\affpsup{2}{h}$, $D_{n+1}\affpsup{2}{h}$) is generalized positive semidefinite whenever $h\geq\frac{3\sqrt{2}}{2}$. This uniform bound is sharp over the whole family, since equality is required when~$n=2$.
\end {Proposition}

\begin{proof}
By Lemma~\ref{1.2}, the individual generalized-semidefinite threshold of an affine matrix~$S\psup{h}$ is the spectral radius of the symmetric nonnegative matrix~$K$ in $\widehat{S\psup{h}}=hI-K$. Equivalently, it is the largest zero of $\det\widehat{S\psup{h}}$.

{\rm (i)} Expanding the affine $B_n^{(1)}$ diagram at its affine leaf gives
\[
 \widehat b_n\affpsup{1}{h}=h\bigl(\widehat b_n\psup{h}-\widehat b_{n-2}\psup{h}\bigr),\qquad n\ge3.
\]
For~$n=3$, this is
$\widehat b_3\affpsup{1}{h}=h^2(h^2-\tfrac{17}{4})$, so~$\lambda_3=\tfrac{\sqrt{17}}{2}$. Since $g$ increases on $[2,\infty)$ and $g(\tfrac{\sqrt{17}}2)=\tfrac{\sqrt{17}}2-2>0$, we have $\epsilon<\tfrac{\sqrt{17}}2$.
If~$n\ge4$ and~$h\ge\tfrac{\sqrt{17}}{2}$, then $h>\epsilon$ and Lemma~\ref{4.10} yields $\widehat b_n\psup{h}>\widehat b_{n-2}\psup{h}>0$. Hence, $\widehat b_n\affpsup{1}{h}>0$ throughout $[\tfrac{\sqrt{17}}{2},\infty)$. Therefore, its largest zero satisfies $\lambda_n<\tfrac{\sqrt{17}}{2}$ for~$n\ge4$, while equality occurs for~$n=3$. This proves the asserted uniform sharp bound.

{\rm (ii)} Let~$K_n$ be the symmetric nonnegative matrix determined by
$\widehat{C_n\affpsup{1}{h}}=hI-K_n$. Its graph is a path on $n+1$ vertices; the two end edges have weight $\tfrac{3}{2}$ and every interior edge has weight $1$. Put~$R\mathrel{:=}\tfrac{3\sqrt2}{2}$ and choose the positive vector
\[
 v=(1,\sqrt2,\ldots,\sqrt2,1)^{\mathsf T}\in\mathbb R^{n+1}.
\]

\Needspace{6\baselineskip}
For~$n=2$, direct multiplication gives
\[
 K_2v=
 \begin{pmatrix}
 0&\tfrac32&0\\
 \tfrac32&0&\tfrac32\\
 0&\tfrac32&0
 \end{pmatrix}
 \begin{pmatrix}1\\\sqrt2\\1\end{pmatrix}
 =\begin{pmatrix}\tfrac{3\sqrt2}{2}\\3\\\tfrac{3\sqrt2}{2}\end{pmatrix}
 =Rv.
\]
Since~$v>0$, the Perron--Frobenius theorem gives $\rho(K_2)=R$. Suppose now that~$n\ge3$.
 At each endpoint, $(K_nv)_i=(\tfrac{3}{2})\sqrt2=Rv_i$. At a vertex adjacent to an endpoint,
$(K_nv)_i=\tfrac{3}{2}+\sqrt2<3=R\sqrt2=Rv_i$, and at every other interior vertex
$(K_nv)_i\le2\sqrt2<3=Rv_i$. Thus, $K_nv\le Rv$ componentwise. The Collatz--Wielandt form of the Perron--Frobenius theorem yields
\[
 \eta_n=\rho(K_n)\le R=\frac{3\sqrt2}{2}
\]
for all $n>1$. The equality already established for~$n=2$ proves that the uniform bound is sharp. By the weighted-graph identification stated above, the paired twisted affine types have permutation-equivalent symmetrizations; hence, the same spectral-radius bounds apply to them. \end{proof}

\begin {Remark} \label{4.17}
For the exact threshold statements in Theorem~\ref{4.4}, Proposition~\ref{4.5'}, and Proposition~\ref{4.5}, replacing the relevant semidefinite condition by its corresponding definite condition replaces the weak threshold inequality by the corresponding strict inequality. Proposition~\ref{4.11} is instead a statement about sharp \emph{uniform sufficient bounds} over affine families and should not be read as giving the individual threshold for every rank.
\end {Remark}

\begin{Lemma}\label{affine-nonsymmetric}
Let~$A$ be an indecomposable symmetrizable affine generalized Cartan matrix. Then, $A$ is generalized positive semidefinite if and only if~$A$ is symmetric.
\end{Lemma}
\begin{proof}
Choose a positive diagonal symmetrizer $D$ such that~$S=DA$ is symmetric. For affine type, $S$ is positive semidefinite of corank one and~$\ker S=\mathbb R u$ for some vector~$u>0$. Since~$Au=0$, one has $u^{\mathsf T}\widehat A u=0$, where~$\widehat A=\tfrac12(A+A^{\mathsf T})$.

If~$\widehat A$ were positive semidefinite, the equality $u^{\mathsf T}\widehat A u=0$ would imply $\widehat A u=0$. But
\[
 2\widehat A u=A^{\mathsf T}u=S D^{-1}u.
\]
Hence, $D^{-1}u\in\ker S=\mathbb Ru$. Since every entry of $u$ is positive, this forces $D$ to be a scalar matrix, and therefore $A=D^{-1}S$ is symmetric. The converse is the standard positive-semidefiniteness of a symmetric affine Cartan matrix. \end{proof}

\begin {Proposition} \label{4.13}
{\rm (1)} An indecomposable generalized Cartan matrix~$H$ is generalized positive semidefinite if and only if, after a simultaneous permutation of rows and columns, it is one of
\[
\begin{gathered}
A_n\psup{2}\ (n\ge1),\quad D_n\psup{2}\ (n\ge4),\quad E_6\psup{2},E_7\psup{2},E_8\psup{2},\\
B_n\psup{2},C_n\psup{2}\ (2\le n\le9),\quad F_4\psup{2},G_2\psup{2},\\
A_1\affpsup{1}{2},\ A_n\affpsup{1}{2}\ (n\ge2),\ D_n\affpsup{1}{2}\ (n\ge4),\ E_6\affpsup{1}{2},E_7\affpsup{1}{2},E_8\affpsup{1}{2}.
\end{gathered}
\]

{\rm (2)} It is generalized positive definite if and only if, after a simultaneous permutation, it is one of
\[
 A_n\psup{2}\ (n\ge1),\quad D_n\psup{2}\ (n\ge4),\quad E_6\psup{2},E_7\psup{2},E_8\psup{2},\quad
 B_n\psup{2},C_n\psup{2}\ (2\le n\le8).
\]
\end {Proposition}

\begin{proof}
The finite nonsymmetric types are controlled by Proposition~\ref{4.8}. Since
\[
 \widehat b_n\psup{2}=2-\frac{n-1}{4},
\]
we have $\widehat b_n\psup{2}>0$ for~$n\le8$, $\widehat b_9\psup{2}=0$, and~$\widehat b_n\psup{2}<0$ for~$n\ge10$. Moreover, $\det\widehat{F_4\psup{2}}=0$ and~$\det\widehat{G_2\psup{2}}=0$. Thus, $B_n\psup{2}$ and~$C_n\psup{2}$ are generalized positive definite for~$2\le n\le8$, while $B_9\psup{2},C_9\psup{2},F_4\psup{2},G_2\psup{2}$ are generalized positive semidefinite but not generalized positive definite.

For the symmetric finite types, the usual Cartan theory gives positive definiteness of
\[
A_n\psup{2},D_n\psup{2},E_6\psup{2},E_7\psup{2},E_8\psup{2}.
\]
The symmetric affine types
\[
 A_1\affpsup{1}{2},\quad
 A_n\affpsup{1}{2},\quad
 D_n\affpsup{1}{2},\quad
 E_6\affpsup{1}{2},\quad
 E_7\affpsup{1}{2},\quad
 E_8\affpsup{1}{2}
\]
are positive semidefinite of corank one. By Lemma~\ref{affine-nonsymmetric}, no nonsymmetric affine type is generalized positive semidefinite.

Finally, an indecomposable generalized Cartan matrix of indefinite type admits a vector~$u>0$ with~$Hu<0$; see \cite[Chapter~4]{Ka90}. Hence, $u^{\mathsf T}Hu<0$, so such a matrix cannot be generalized positive semidefinite. These three cases exhaust the finite, affine, and indefinite possibilities and prove both assertions.
\end{proof}

\section[Determinant sequences associated with h-Shuhan matrices]{Determinant sequences associated with~$h$-Shuhan matrices}\label{s4}
The classification in Section~\ref{s3} concerns $0\le h<2$, but the determinant recurrences themselves are defined for every $h$. This section develops them as a second, complementary theme. The congruence results require integral $h$, whereas the asymptotic results in the later subsections assume $h>2$, where the characteristic roots $\alpha_h,\beta_h$ are real and satisfy $\alpha_h>1>\beta_h>0$.

\subsection{Congruence periodicity of the determinant sequences}
In this subsection, $h\in\mathbb Z$ and~$\Delta\mathrel{:=}h^2-4$. The determinant identities used here are polynomial identities in $h$, so the formulas established for~$h\ge0$ extend to every real $h$. For an odd prime $p$, write
\begin{equation}\label{e5.1}
\left(\frac{\Delta}{p}\right):=
\begin{cases}
 1,& \Delta\text{ is a nonzero quadratic residue modulo }p,\\
-1,& \Delta\text{ is a quadratic nonresidue modulo }p,\\
 0,& p\mid\Delta.
\end{cases}
\end{equation}

\begin{Proposition}\label{9.01}
Let~$p$ be an odd prime.
\begin{parts}
\item $a_{p-1}\psup{h}\equiv\left(\frac{\Delta}{p}\right)\pmod p$ and~$b_p\psup{h}\equiv h\pmod p$.
\item If~$\left(\frac{\Delta}{p}\right)=1$, then $a_{p-1+i}\psup{h}\equiv a_i\psup{h}\pmod p$ for every $i\ge0$.
\item If~$\left(\frac{\Delta}{p}\right)=-1$, then $a_{p+1+i}\psup{h}\equiv a_i\psup{h}\pmod p$ for every $i\ge0$.
\item If~$\left(\frac{\Delta}{p}\right)=0$, then the $a$-sequence has period $p$ when~$h\equiv2\pmod p$ and period $2p$ when~$h\equiv-2\pmod p$.
\end{parts}
\end{Proposition}
\begin{proof} {\rm (i)} By Proposition~\ref{4.001},
\[
 2^{p-1}a_{p-1}\psup{h}
 =\sum_{0\le k<\tfrac p2}\binom{p}{2k+1}h^{p-1-2k}\Delta^k.
\]
Modulo $p$, every intermediate binomial coefficient $\binom pk$ with~$1\le k\le p-1$ vanishes. Hence, using Fermat's little theorem $2^{p-1}\equiv1\pmod p$,
\[
 a_{p-1}\psup{h}\equiv \Delta^{\tfrac{p-1}{2}}\equiv\left(\frac{\Delta}{p}\right)\pmod p
\]
by Euler's criterion. Similarly,
\[
 2^{p-1}b_p\psup{h}
 =\sum_{0\le k\le\tfrac p2}\binom{p}{2k}h^{p-2k}\Delta^k
 \equiv h^p\equiv h\pmod p,
\]
which proves (i).

{\rm (ii)} Assume $\left(\frac{\Delta}{p}\right)=1$. Part (i) gives $a_{p-1}\psup{h}\equiv1=a_0\psup{h}$. From
\[
 b_p\psup{h}=a_p\psup{h}-a_{p-2}\psup{h}=2a_p\psup{h}-ha_{p-1}\psup{h}
\]
we obtain
\[
 2a_p\psup{h}=b_p\psup{h}+ha_{p-1}\psup{h}\equiv h+h=2h\pmod p,
\]
so~$a_p\psup{h}\equiv h=a_1\psup{h}$. Thus, the two consecutive pairs $(a_{p-1}\psup{h},a_p\psup{h})$ and~$(a_0\psup{h},a_1\psup{h})$ are congruent modulo $p$. Since~$a_{n+1}\psup{h}=ha_n\psup{h}-a_{n-1}\psup{h}$, induction now yields
$a_{p-1+i}\psup{h}\equiv a_i\psup{h}\pmod p$ for every $i\ge0$.

{\rm (iii)} Assume $\left(\frac{\Delta}{p}\right)=-1$. The same identity gives
$2a_p\psup{h}\equiv h-h=0\pmod p$, hence~$a_p\psup{h}\equiv0$. Therefore, \[
 a_{p+1}\psup{h}\equiv-a_{p-1}\psup{h}\equiv1=a_0\psup{h},
 \qquad
 a_{p+2}\psup{h}\equiv h=a_1\psup{h}.
\]
The recurrence applied to these two consecutive initial congruences gives
$a_{p+1+i}\psup{h}\equiv a_i\psup{h}\pmod p$ for every $i\ge0$.

{\rm (iv)} If~$\left(\frac{\Delta}{p}\right)=0$, then $h\equiv2$ or $h\equiv-2\pmod p$.

\noindent If~$h\equiv2\pmod p$, the recurrence with~$(a_0\psup{2},a_1\psup{2})=(1,2)$ gives $a_n\psup{h}\equiv a_n\psup{2}=n+1\pmod p$. Hence, \[
 (a_p\psup{h},a_{p+1}\psup{h})\equiv(1,2)\equiv(a_0\psup{h},a_1\psup{h})\pmod p,
\]
and the period $p$ follows. If~$h\equiv-2\pmod p$, then
$a_n\psup{h}\equiv(-1)^n(n+1)\pmod p$, so
\[
 (a_{2p}\psup{h},a_{2p+1}\psup{h})\equiv(1,-2)\equiv(a_0\psup{h},a_1\psup{h})\pmod p.
\]
The recurrence therefore gives the period $2p$. This proves (iv). \end{proof}

\begin{Proposition}\label{9.02}
Let~$p$ be an odd prime.
\begin{parts}
\item If~$\left(\frac{\Delta}{p}\right)=1$, then $b_{p-1+i}\psup{h}\equiv b_i\psup{h}\pmod p$ for every $i\ge1$.
\item If~$\left(\frac{\Delta}{p}\right)=-1$, then $b_{p+1+i}\psup{h}\equiv b_i\psup{h}\pmod p$ for every $i\ge1$.
\item If~$\left(\frac{\Delta}{p}\right)=0$, then the $b$-sequence has period $1$ when~$h\equiv2\pmod p$ and period $2$ when~$h\equiv-2\pmod p$.
\end{parts}
\end{Proposition}
\begin{proof}
The recurrence $b_{n+1}\psup{h}=hb_n\psup{h}-b_{n-1}\psup{h}$ holds for~$n\ge2$. We shall repeatedly use the identity
\[
 2b_{p+1}\psup{h}=\Delta a_{p-1}\psup{h}+h b_p\psup{h},
\]
which follows from $b_n\psup{h}=ha_{n-1}\psup{h}-2a_{n-2}\psup{h}$ and the recurrence for~$a_n\psup{h}$.

{\rm (i)} If~$\left(\frac{\Delta}{p}\right)=1$, Proposition~\ref{9.01}(i) gives
$b_p\psup{h}\equiv h=b_1\psup{h}$ and~$a_{p-1}\psup{h}\equiv1$. Hence, this identity yields
\[
 2b_{p+1}\psup{h}\equiv\Delta+h^2=2(h^2-2)=2b_2\psup{h}\pmod p.
\]
Since~$p$ is odd, $b_{p+1}\psup{h}\equiv b_2\psup{h}$. The two consecutive initial congruences now imply
$b_{p-1+i}\psup{h}\equiv b_i\psup{h}\pmod p$ for all $i\ge1$.

{\rm (ii)} If~$\left(\frac{\Delta}{p}\right)=-1$, then $a_{p-1}\psup{h}\equiv-1$ and~$b_p\psup{h}\equiv h$. The same identity gives
$2b_{p+1}\psup{h}\equiv-\Delta+h^2=4$, hence~$b_{p+1}\psup{h}\equiv2$. Consequently, \[
 b_{p+2}\psup{h}\equiv2h-h=h=b_1\psup{h},
 \qquad
 b_{p+3}\psup{h}\equiv h^2-2=b_2\psup{h}.
\]
The recurrence gives $b_{p+1+i}\psup{h}\equiv b_i\psup{h}\pmod p$ for all $i\ge1$.

{\rm (iii)} If~$\left(\frac{\Delta}{p}\right)=0$, then $h\equiv\pm2\pmod p$. For~$h\equiv2$, the recurrence gives $b_n\psup{h}\equiv2$ for all $n$, hence period $1$. For~$h\equiv-2$, it gives $b_n\psup{h}\equiv2(-1)^n$, hence period $2$. This is exactly the assertion. \end{proof}

\begin {Proposition} \label{9.05}
Let~$p$ be an odd prime.

{\rm (i)} If~$\left(\frac{\Delta}{p}\right)=1$, then
$d_i\psup{h}\equiv d_{p-1+i}\psup{h}\pmod p$ for every $i\ge3$.

{\rm (ii)} If~$\left(\frac{\Delta}{p}\right)=-1$, then
$d_i\psup{h}\equiv d_{p+1+i}\psup{h}\pmod p$ for every $i\ge3$.

{\rm (iii)} If~$\left(\frac{\Delta}{p}\right)=0$, then
$d_i\psup{h}\equiv d_{1+i}\psup{h}$ or $d_i\psup{h}\equiv d_{2+i}\psup{h}\pmod p$ for every $i\ge3$, according as $h\equiv2$ or $-2\pmod p$.
\end {Proposition}
\begin{proof}
By \eqref{e3}, $d_i\psup{h}=h b_{i-1}\psup{h}$. All three parts follow immediately from Proposition~\ref{9.02} after multiplying its congruences by~$h$. \end{proof}

\begin {Proposition} \label{9.03}
{\rm (i)} If~$2\mid\Delta$ (equivalently, $h$ is even), then
\[
 a_i\psup{h}\equiv a_{i+2}\psup{h}\pmod2,\qquad
 b_{i+1}\psup{h}\equiv b_{i+2}\psup{h}\pmod2
\]
for every $i\ge0$.

{\rm (ii)} If~$2\nmid\Delta$ (equivalently, $h$ is odd), then
\[
 a_i\psup{h}\equiv a_{i+3}\psup{h}\pmod2,\qquad
 b_{i+1}\psup{h}\equiv b_{i+4}\psup{h}\pmod2,
\]
and
\[
 d_{i+3}\psup{h}\equiv d_{i+6}\psup{h}\pmod2
\]
for every $i\ge0$ for which the displayed terms are defined.
\end {Proposition}
\begin{proof}
If~$h$ is even, the recurrence gives $a_n\psup{h}\equiv a_{n-2}\psup{h}\pmod2$. Also $b_1\psup{h}=h\equiv0$ and~$b_2\psup{h}=h^2-2\equiv0\pmod2$, so the recurrence forces $b_n\psup{h}\equiv0\pmod2$ for every $n$, proving (i).

If~$h$ is odd, both $a_n\psup{h}$ and~$b_n\psup{h}$ satisfy $x_n\equiv x_{n-1}+x_{n-2}\pmod2$. Their initial residues are respectively
$(a_0\psup{h},a_1\psup{h})\equiv(1,1)$ and~$(b_1\psup{h},b_2\psup{h})\equiv(1,1)$, so both residue sequences repeat the pattern $1,1,0$ with period $3$. Finally, $d_n\psup{h}=h b_{n-1}\psup{h}\equiv b_{n-1}\psup{h}\pmod2$, and the same period-$3$ conclusion follows for~$d_n\psup{h}$. \end{proof}

\begin {Remark} \label{9.04}
Propositions~\ref{9.01}, \ref{9.02}, \ref{9.05}, and \ref{9.03} give explicit periods modulo every prime. We do not claim that the displayed periods are always minimal.
\end {Remark}
\subsection[Determinant sequences for h > 2]{The sequences $a_{n}\psup{h},b_{n}\psup{h}$ and~$d_{n}\psup{h}$ with~$h>2$}
If~$h>2$, then $0<\beta_{h}<1$.

\begin{Proposition}\label{4.02}
Assume $h>2$. Then, the following statements hold.
\begin{parts}
\item For each of the three determinant sequences $x_n\psup{h}=a_n\psup{h}$, $b_n\psup{h}$, or $d_n\psup{h}$, we have, whenever the fractions are defined,
\[
 \frac{x_{n-1}\psup{h}}{x_n\psup{h}}\longrightarrow\beta_h\qquad(n\to\infty).
\]
\item Each of the three determinant sequences $a_n\psup{h}$, $b_n\psup{h}$, and $d_n\psup{h}$ is strictly increasing for the indices where it is defined. Moreover, the ratios $\tfrac{a_{n-1}\psup{h}}{a_n\psup{h}}$ are strictly increasing.
\item The ratios $\tfrac{b_{n-1}\psup{h}}{b_n\psup{h}}$ and $\tfrac{d_{n-1}\psup{h}}{d_n\psup{h}}$ form strictly decreasing sequences whenever the fractions are defined.
\end{parts}
\end{Proposition}
\begin{proof}
Using the explicit formulas in Proposition~\ref{4.01} and $\alpha_h\beta_h=1$,
\[
 \frac{a_{n-1}\psup{h}}{a_n\psup{h}}
 =\frac{\beta_h-\beta_h^{2n+1}}{1-\beta_h^{2n+2}},\qquad
 \frac{b_{n-1}\psup{h}}{b_n\psup{h}}
 =\frac{\beta_h+\beta_h^{2n-1}}{1+\beta_h^{2n}},
\]
and
\[
 \frac{d_{n-1}\psup{h}}{d_n\psup{h}}
 =\frac{\beta_h+\beta_h^{2n-3}}{1+\beta_h^{2n-2}}.
\]
Because $0<\beta_h<1$, all three ratios tend to $\beta_h$, proving (i).  The displayed formulas also show that all three sequences are positive and that each consecutive-term ratio is $<1$, hence~$a_n\psup{h},b_n\psup{h},d_n\psup{h}$ are strictly increasing.

It remains to determine the monotonicity of the ratios.  Since~$f(a\psup{h})=-1$,
\[
 a_n\psup{h} a_{n-2}\psup{h}=(a_{n-1}\psup{h})^2-1<(a_{n-1}\psup{h})^2,
\]
so
\[
 \frac{a_{n-2}\psup{h}}{a_{n-1}\psup{h}}<\frac{a_{n-1}\psup{h}}{a_n\psup{h}};
\]
thus $\tfrac{a_{n-1}\psup{h}}{a_n\psup{h}}$ increases.  On the other hand, $f(b\psup{h})=(\alpha_h-\beta_h)^2>0$ gives
\[
 b_n\psup{h} b_{n-2}\psup{h}=(b_{n-1}\psup{h})^2+f(b\psup{h})>(b_{n-1}\psup{h})^2,
\]
whence $\tfrac{b_{n-1}\psup{h}}{b_n\psup{h}}<\tfrac{b_{n-2}\psup{h}}{b_{n-1}\psup{h}}$.  Finally,
$d_n\psup{h}=hb_{n-1}\psup{h}$, so the same conclusion holds for the $d$-ratios.  This proves (ii) and (iii). \end{proof}

\begin{Proposition}\label{4.03}
Assume $h>2$ and~$n\ge2$. Then, we have
\begin{parts}
\item
\[
 a_n\affpsup{1}{h}=h a_n\psup{h}-2a_{n-1}\psup{h}-2
 =\alpha_h^{n+1}+\beta_h^{n+1}-2.
\]
\item
\[
 \frac{a_{n-1}\affpsup{1}{h}}{a_n\affpsup{1}{h}}\longrightarrow\beta_h\qquad(n\to\infty).
\]
\end{parts}
\end{Proposition}
\begin{proof}
Expanding the determinant of the cyclic matrix~$A_n\affpsup{1}{h}$ along the affine row (or column) gives
\[
 a_n\affpsup{1}{h}=h a_n\psup{h}-2a_{n-1}\psup{h}-2.
\]
Using $b_{n+1}\psup{h}=h a_n\psup{h}-2a_{n-1}\psup{h}$ from \eqref{e3}, we obtain
\[
 a_n\affpsup{1}{h}=b_{n+1}\psup{h}-2
 =\alpha_h^{n+1}+\beta_h^{n+1}-2,
\]
which proves (i).  For (ii), divide numerator and denominator by~$\alpha_h^{n+1}$ and use $\alpha_h^{-1}=\beta_h$ and~$0<\beta_h<1$:
\[
 \frac{a_{n-1}\affpsup{1}{h}}{a_n\affpsup{1}{h}}
 =\frac{\alpha_h^n+\beta_h^n-2}{\alpha_h^{n+1}+\beta_h^{n+1}-2}
 \longrightarrow\frac1{\alpha_h}=\beta_h\qquad(n\to\infty).
\]
\end{proof}

\subsection{Two-index determinant sequences}

For~$1\le i\le n-1$, let~$A_{i,n}\psup{h}$ be obtained from $A_n\psup{h}$ by replacing both entries in positions $(i,i+1)$ and~$(i+1,i)$ by~$-2$, and let~$B_{i,n}\psup{h}$ be obtained from $A_n\psup{h}$ by replacing only the entry in position $(i,i+1)$ by~$-2$.

Set
\[
 a_{i,n}\psup{h}:=\det A_{i,n}\psup{h},
 \qquad
 b_{i,n}\psup{h}:=\det B_{i,n}\psup{h}.
\]
We also put~$a_{i,i}\psup{h}:=a_i\psup{h}$ and~$b_{i,i}\psup{h}:=a_i\psup{h}$. Expansion at the last row gives the second-order recurrences
\begin{equation}\label{e:two-index-rec}
\begin{aligned}
 a_{i,i+1}\psup{h}&=h a_i\psup{h}-4a_{i-1}\psup{h},
 & a_{i,n}\psup{h}&=h a_{i,n-1}\psup{h}-a_{i,n-2}\psup{h},\\
 b_{i,i+1}\psup{h}&=h a_i\psup{h}-2a_{i-1}\psup{h},
 & b_{i,n}\psup{h}&=h b_{i,n-1}\psup{h}-b_{i,n-2}\psup{h},
\end{aligned}
\end{equation}
where the recurrences in the second column apply for~$1\le i\le n-2$.
Reflection of the path gives the first symmetry below; for~$B_{i,n}\psup{h}$ one additionally uses invariance of the determinant under transpose, which reverses the modified edge. Hence,
\[
 a_{i,n}\psup{h}=a_{n-i,n}\psup{h},
 \qquad
 b_{i,n}\psup{h}=b_{n-i,n}\psup{h}.
\]

For each fixed $i$, the sequences $(a_{i,n}\psup{h})_{n\ge i}$ and
$(b_{i,n}\psup{h})_{n\ge i}$ satisfy the recurrence
$y_n=hy_{n-1}-y_{n-2}$ for $n\ge i+2$. We write
\[
 F_i^a(h):=f\bigl((a_{i,n}\psup{h})_{n\ge i}\bigr),\qquad
 F_i^b(h):=f\bigl((b_{i,n}\psup{h})_{n\ge i}\bigr),
\]
where $f(y)=h y_ny_{n-1}-y_n^2-y_{n-1}^2$ is the
recurrence invariant of the whole sequence $y=(y_n)$.
Thus $F_i^a(h)$ and $F_i^b(h)$ do not depend on $n$.
The next results determine positivity and asymptotic behavior
in the ranges needed later.

\begin {Lemma} \label{5.1}
Let~$h>2$. Then, \[
 \lim_{i\to\infty}
 \frac{F_{i}^a(h)}{F_{i+1}^a(h)}
 =
 \begin{cases}
   \beta_h^2,& h\ne\frac52,\\
   1,& h=\frac52,
 \end{cases}
\]
and
\[
 \lim_{i\to\infty}
 \frac{F_{i}^b(h)}{F_{i+1}^b(h)}
 =
 \begin{cases}
   \beta_h^2,& h\ne\frac{3\sqrt2}{2},\\
   1,& h=\frac{3\sqrt2}{2}.
 \end{cases}
\]
More generally, for every admissible $n\ge i+2$, the recurrence invariant satisfies
\[
\begin{aligned}
 F_{i}^a(h)&=4h a_i\psup{h} a_{i-1}\psup{h}-(a_i\psup{h})^2-16(a_{i-1}\psup{h})^2,\\
 F_{i}^b(h)&=2h a_i\psup{h} a_{i-1}\psup{h}-(a_i\psup{h})^2-4(a_{i-1}\psup{h})^2.
\end{aligned}
\]
\end {Lemma}
\begin{proof}
Using the recurrence in \eqref{e:two-index-rec}, the quantity $f$ is independent of the second index once the first index $i$ is fixed. Evaluating it at~$n=i+2$ gives the two explicit formulas displayed in the statement.
By Proposition~\ref{4.02}, $\tfrac{a_{i-1}\psup{h}}{a_i\psup{h}}\to\beta_h$ and~$\tfrac{a_i\psup{h}}{a_{i+1}\psup{h}}\to\beta_h$ as $i\to\infty$. If
\[
 P_a(\beta_h):=4h\beta_h-1-16\beta_h^2=3(1-4\beta_h^2)\ne0,
\]
then dividing the first consecutive ratio by~$(a_{i+1}\psup{h})^2$ gives
\[
 \frac{F_{i}^a(h)}{F_{i+1}^a(h)}
 =\frac{4h\frac{a_{i-1}\psup{h}}{a_i\psup{h}}-1-16(\frac{a_{i-1}\psup{h}}{a_i\psup{h}})^2}
 {4h\frac{a_i\psup{h}}{a_{i+1}\psup{h}}-1-16(\frac{a_i\psup{h}}{a_{i+1}\psup{h}})^2}
 \left(\frac{a_i\psup{h}}{a_{i+1}\psup{h}}\right)^2
 \longrightarrow\beta_h^2\qquad(i\to\infty).
\]
The exceptional equation $P_a(\beta_h)=0$ is equivalent to $\beta_h=\tfrac{1}{2}$, hence~$h=\alpha_h+\beta_h=\tfrac{5}{2}$. In this case the closed form $a_i\psup{\tfrac{5}{2}}=\frac{4^{i+1}-1}{3\,2^i}$ gives
\[
 F_{i}^a(\tfrac{5}{2})=4-\frac{5}{4^i}\longrightarrow4\qquad(i\to\infty),
\]
so the consecutive ratio tends to $1$.

For the $b$-sequence, the same argument uses
$P_b(\beta_h)\mathrel{:=}2h\beta_h-1-4\beta_h^2=1-2\beta_h^2$. Thus, the generic limit is $\beta_h^2$, while the exceptional case is $\beta_h=\tfrac{1}{\sqrt2}$, equivalently $h=\tfrac{3\sqrt2}{2}$. Using $\alpha_h=\sqrt2$, $\beta_h=\tfrac{1}{\sqrt2}$ in the closed form for~$a_i\psup{h}$ gives explicitly
\[
 F_{i}^b(\tfrac{3\sqrt2}{2})=2-\frac{3}{2^i}\longrightarrow2\qquad(i\to\infty),
\]
so the consecutive ratio again tends to $1$. \end{proof}

\begin{Lemma}\label{5.2}
\begin{parts}
\item If~$h>2$, then
\[
 \frac{\beta_h}{2}\le\frac{a_{i-1}\psup{h}}{a_i\psup{h}}<\beta_h
 \qquad(i\ge1).
\]
\item If~$h\ge\tfrac52$, then $\beta_h\le\tfrac{\alpha_h}{4}$ and~$F_{i}^a(h)>0$ for~$1\le i\le n-2$.
\item If~$h\ge\tfrac{3\sqrt2}{2}$, then
\[
 \frac{2\beta_h}{3}\le\frac{a_{i-1}\psup{h}}{a_i\psup{h}}<\beta_h
 \qquad(i\ge1),
\]
and equality on the left occurs only for~$h=\tfrac{3\sqrt2}{2}$ and~$i=1$.
\item If~$h\ge\tfrac{3\sqrt2}{2}$, then $\beta_h\le\tfrac{\alpha_h}{2}$ and~$F_{i}^b(h)>0$ for~$1\le i\le n-2$.
\end{parts}
\end{Lemma}
\begin{proof}
{\rm (i)} Proposition~\ref{4.02} gives
\[
 r_i:=\frac{a_{i-1}\psup{h}}{a_i\psup{h}}
 =\beta_h\frac{1-\beta_h^{2i}}{1-\beta_h^{2i+2}}<\beta_h.
\]
For the lower bound, $r_i\ge\tfrac{\beta_h}{2}$ is equivalent to
$1\ge(2-\beta_h^2)\beta_h^{2i}$. Since~$i\ge1$ and~$0<\beta_h<1$,
\[
 (2-\beta_h^2)\beta_h^{2i}
 \le (2-\beta_h^2)\beta_h^2
 =1-(1-\beta_h^2)^2\le1,
\]
which proves (i).

{\rm (ii)} If~$h\ge\tfrac52$, then $\alpha_h\ge\alpha_{\tfrac{5}{2}}=2$. Hence, $\beta_h=\alpha_h^{-1}\le\tfrac{\alpha_h}{4}$. Moreover,
\[
 (1-\beta_h^{2i+2})\frac{\alpha_h}{4}
 =\frac{\alpha_h}{4}-\frac{\beta_h^{2i+1}}{4}
 >\frac{\alpha_h}{4}-\beta_h^{2i+1}
 \ge\beta_h-\beta_h^{2i+1},
\]
so~$r_i<\tfrac{\alpha_h}{4}$. Part (i) gives $r_i\ge\tfrac{\beta_h}{2}>\tfrac{\beta_h}{4}$. Since the roots of
$4hx-1-16x^2=0$ are $\tfrac{\beta_h}{4}$ and~$\tfrac{\alpha_h}{4}$, it follows that
\[
 \frac{F_{i}^a(h)}{(a_i\psup{h})^2}=4hr_i-1-16r_i^2>0,
\]
proving (ii).

{\rm (iii)} If~$h\geq\frac{3\sqrt{2}}{2}$, then $\sqrt{h^{2}-4}\geq \tfrac h3$, and therefore
\[
 \alpha_h=\frac{h+\sqrt{h^2-4}}2\ge\frac{2h}{3},
 \qquad
 \beta_h=\alpha_h^{-1}\le\frac{3}{2h}.
\]
For~$i=1$, we have

\[
 \frac{a_0\psup{h}}{a_1\psup{h}}=\frac1h\ge\frac{2\beta_h}{3}.
\]
  Equality holds here precisely when~$h=\tfrac{3\sqrt2}{2}$.  If~$i\ge2$, then the recurrence $ha_{i-1}\psup{h}=a_i\psup{h}+a_{i-2}\psup{h}$ and positivity of $a_{i-2}\psup{h}$ give
\[
 \frac{a_{i-1}\psup{h}}{a_i\psup{h}}
 =\frac{a_i\psup{h}+a_{i-2}\psup{h}}{h a_i\psup{h}}
 >\frac1h\ge\frac{2\beta_h}{3}.
\]
The upper bound $r_i<\beta_h$ was established in part~(i). This proves part~(iii), including its equality statement.

{\rm (iv)} Since~$h\ge\tfrac{3\sqrt2}{2}$, one has $\beta_h\le\tfrac{1}{\sqrt2}$, and hence
$\beta_h\le\tfrac{\alpha_h}{2}$. By part {\rm (iii)},
\[
 \frac{\beta_h}{2}<\frac{2\beta_h}{3}\le\frac{a_{i-1}\psup{h}}{a_i\psup{h}}
 <\beta_h\le\frac{\alpha_h}{2}.
\]
The two roots of $2hx-1-4x^2=0$ are $\tfrac{\beta_h}{2}$ and~$\tfrac{\alpha_h}{2}$. Therefore, \[
 \frac{F_{i}^b(h)}{(a_i\psup{h})^2}
 =2h\frac{a_{i-1}\psup{h}}{a_i\psup{h}}-1-4\left(\frac{a_{i-1}\psup{h}}{a_i\psup{h}}\right)^2>0,
\]
which proves the assertion. \end{proof}

Next, define $\gamma_1$ and~$\gamma_2$ as follows.  Let~$\gamma_1\approx2.561924625$ be the unique real root $>2$ (equivalently, the largest real root) of
\[
 3h^4-6h^3-15h^2+34h-17=0,
\]
and let~$\gamma_2\approx2.173868929$ be the unique real root $>2$ of
\[
 h^4-2h^3-3h^2+8h-5=0.
\]

The asserted uniqueness of these roots is elementary. The polynomial
\[
 p_1(h):=3h^4-6h^3-15h^2+34h-17
\]
satisfies $p_1(2)=-9$, $p_1''(h)=36h^2-36h-30>0$ for~$h\ge2$, and~$p_1'(2)=-2$. Hence, $p_1'$ has at most one zero on $[2,\infty)$; $p_1$ first decreases and then increases to $+\infty$, so it has exactly one zero greater than $2$. Likewise, for
\[
 p_2(h):=h^4-2h^3-3h^2+8h-5,
\]
we have $p_2(2)=-1$, $p_2'(2)=4$, and~$p_2''(h)=12h^2-12h-6>0$ for~$h\ge2$, so~$p_2$ is strictly increasing on $[2,\infty)$ and has exactly one zero there.

Since
\[
 p_1(\tfrac52)=-\tfrac{37}{16}<0,\qquad p_2(\tfrac{3\sqrt2}2)=\tfrac74-\tfrac{3\sqrt2}2<0,
\]
uniqueness gives $\gamma_1>\tfrac52$ and $\gamma_2>\tfrac{3\sqrt2}2$, as required in Lemma~\ref{5.2}(ii),(iv).

For $h>2$, put $q(h):=h^2-h-1$. Using $\alpha_h+\beta_h=h$ and $\alpha_h\beta_h=1$ gives
\[
\begin{aligned}
 p_1(h)&=-[q(h)-4(h-1)\alpha_h][q(h)-4(h-1)\beta_h],\\
 p_2(h)&=-[q(h)-2(h-1)\alpha_h][q(h)-2(h-1)\beta_h].
\end{aligned}
\]
Since $\alpha_h>\tfrac h2$, the factors containing $\alpha_h$ are negative. Hence, $q(h)-4(h-1)\beta_h$ and $q(h)-2(h-1)\beta_h$ vanish at $\gamma_1,\gamma_2$, respectively, and are strictly positive for $h>\gamma_1$ and $h>\gamma_2$, respectively. These are the inequalities used below.

\begin{Proposition}\label{5.3}
\begin{parts}
\item If~$h>\gamma_1$, then $a_{i,n}\psup{h}>0$ for every admissible pair $1\le i\le n-1$, and for each fixed $i$ the sequence $n\mapsto a_{i,n}\psup{h}$ is strictly increasing for~$n>i$.
\item If~$h>\gamma_2$, then $b_{i,n}\psup{h}>0$ for every admissible pair $1\le i\le n-1$, and for each fixed $i$ the sequence $n\mapsto b_{i,n}\psup{h}$ is strictly increasing for~$n>i$.
\end{parts}
\end{Proposition}
\begin{proof}
Fix $i$ and write $r_i=\tfrac{a_{i-1}\psup{h}}{a_i\psup{h}}$. By Lemma~\ref{5.2}(i), $r_i<\beta_h$.

For part~(i),
\[
 \frac{a_{i,i+1}\psup{h}}{a_i\psup{h}}=h-4r_i>h-4\beta_h=2\sqrt{h^2-4}-h>0
\]
when~$h>\gamma_1$. Moreover,
\[
 \frac{a_{i,i+2}\psup{h}-a_{i,i+1}\psup{h}}{a_i\psup{h}}
 =(h^2-h-1)-4(h-1)r_i
 >(h^2-h-1)-4(h-1)\beta_h>0,
\]
where the last inequality is the defining property of $\gamma_1$. Thus, $a_{i,i+2}\psup{h}>a_{i,i+1}\psup{h}>0$. If
$D_n:=a_{i,n}\psup{h}-a_{i,n-1}\psup{h}$, then the recurrence gives
\[
 D_n=(h-2)a_{i,n-1}\psup{h}+D_{n-1}.
\]
Since~$h>2$, induction shows simultaneously that every $a_{i,n}\psup{h}$ and every $D_n$ is positive. Hence, $a_{i,n}\psup{h}$ is strictly increasing in $n$.

For part~(ii),
\[
 \frac{b_{i,i+1}\psup{h}}{a_i\psup{h}}=h-2r_i>h-2\beta_h=\sqrt{h^2-4}>0.
\]
If~$h>\gamma_2$, then
\[
 \frac{b_{i,i+2}\psup{h}-b_{i,i+1}\psup{h}}{a_i\psup{h}}
 =(h^2-h-1)-2(h-1)r_i
 >(h^2-h-1)-2(h-1)\beta_h>0.
\]
The same difference recurrence,
\[
 b_{i,n}\psup{h}-b_{i,n-1}\psup{h}
 =(h-2)b_{i,n-1}\psup{h}+\bigl(b_{i,n-1}\psup{h}-b_{i,n-2}\psup{h}\bigr),
\]
then proves positivity and strict increase by induction.
\end{proof}

\begin{Proposition}\label{5.4}
For each fixed $i\ge1$:
\begin{parts}
\item if~$h>\gamma_1$, then
\[
 \frac{a_{i,n-1}\psup{h}}{a_{i,n}\psup{h}}\longrightarrow\beta_h\qquad(n\to\infty);
\]
\item if~$h>\gamma_2$, then
\[
 \frac{b_{i,n-1}\psup{h}}{b_{i,n}\psup{h}}\longrightarrow\beta_h\qquad(n\to\infty).
\]
\end{parts}
\end{Proposition}
\begin{proof}

Fix $i\ge1$. In case~(i), let $y_n=a_{i,n}\psup{h}$ and assume $h>\gamma_1$; in case~(ii), let $y_n=b_{i,n}\psup{h}$ and assume $h>\gamma_2$. Lemma~\ref{5.2}(ii),(iv) and Proposition~\ref{5.3} give $f(y)>0$ and $0<y_{n-1}<y_n$ for $n\ge i+2$. The recurrence yields
\[
 f(y)=y_ny_{n-2}-y_{n-1}^2>0,\qquad
 r_n^{-1}=h-r_{n-1}<h,
\]
where $r_n:=\tfrac{y_{n-1}}{y_n}$ and $n\ge i+3$. Thus, $r_n$ decreases in $(0,1)$ and has a limit $L\ge\tfrac1h>0$. Passing to the limit gives $L^{-1}=h-L$. The roots of this equation are $\alpha_h,\beta_h$; since $L\le1<\alpha_h$, we obtain $L=\beta_h$ in both cases.

\end{proof}

Define
\[
 w_n:=\begin{cases}1,& n\text{ even},\\16,& n\text{ odd},\end{cases}
 \qquad
 w_n':=\begin{cases}1,& n\text{ even},\\4,& n\text{ odd}.
 \end{cases}
\]
Let~$u_0:=h$ and~$u_0':=h$. For~$n\ge1$, set
\[
 u_n:=\begin{cases}
 a_{n+1}\psup{h}-4a_{n-1}\psup{h},& n\text{ even},\\
 a_n\psup{h},& n\text{ odd},
 \end{cases}
 \qquad
 u_n':=\begin{cases}
 a_{n+1}\psup{h}-2a_{n-1}\psup{h},& n\text{ even},\\
 a_n\psup{h},& n\text{ odd}.
 \end{cases}
\]

\begin{Theorem}\label{5.6}
Assume $h>2$ and~$1\le i<n-1$. Then, \begin{parts}
\item $F_{i}^a(h)=3u_{i-1}u_i-w_i$;
\item $F_{i}^b(h)=u_{i-1}'u_i'-w_i'$.
\end{parts}
\end{Theorem}
\begin{proof}
By Lemma~\ref{5.1},
\begin{align}
 F_{i}^a(h)
 &=4h a_i\psup{h} a_{i-1}\psup{h}-(a_i\psup{h})^2-16(a_{i-1}\psup{h})^2,\label{eq:Fai}\\
 F_{i}^b(h)
 &=2h a_i\psup{h} a_{i-1}\psup{h}-(a_i\psup{h})^2-4(a_{i-1}\psup{h})^2.\label{eq:Fbi}
\end{align}
We compare these expressions with the definitions of $u_i,u_i'$ and use only the identity $f(a\psup{h})=-1$ from Proposition~\ref{4.01}(i). For $i=1$, the same identity holds at index zero because $a_{-1}\psup{h}=0$ and $a_0\psup{h}=1$.

For part~(i), suppose first that~$i$ is odd. Then, $u_{i-1}=a_i\psup{h}-4a_{i-2}\psup{h}$, $u_i=a_i\psup{h}$, and~$w_i=16$. Subtracting $3u_{i-1}u_i-w_i$ from \eqref{eq:Fai}, and using $a_i\psup{h}=h a_{i-1}\psup{h}-a_{i-2}\psup{h}$, gives
\[
 F_{i}^a(h)-\bigl(3u_{i-1}u_i-w_i\bigr)
 =16\bigl(f(a\psup{h})+1\bigr)=0.
\]
If~$i$ is even, then $u_{i-1}=a_{i-1}\psup{h}$, $u_i=a_{i+1}\psup{h}-4a_{i-1}\psup{h}$, and~$w_i=1$, and the same calculation gives
\[
 F_{i}^a(h)-\bigl(3u_{i-1}u_i-w_i\bigr)
 =f(a\psup{h})+1=0.
\]
This proves (i).

For part~(ii), if~$i$ is odd, then
$u_{i-1}'=a_i\psup{h}-2a_{i-2}\psup{h}$, $u_i'=a_i\psup{h}$, and~$w_i'=4$. Equation~\eqref{eq:Fbi} yields
\[
 F_{i}^b(h)-\bigl(u_{i-1}'u_i'-w_i'\bigr)
 =4\bigl(f(a\psup{h})+1\bigr)=0.
\]
If~$i$ is even, then $u_{i-1}'=a_{i-1}\psup{h}$, $u_i'=a_{i+1}\psup{h}-2a_{i-1}\psup{h}$, and~$w_i'=1$, so
\[
 F_{i}^b(h)-\bigl(u_{i-1}'u_i'-w_i'\bigr)
 =f(a\psup{h})+1=0.
\]
Hence, (ii) follows as well.
\end{proof}

\begin {Remark} \label{5.6'}
Regard $a_n\psup{h}$ as a polynomial in the indeterminate $h$. If~$n$ is odd, its recurrence shows inductively that~$h$ divides $a_n\psup{h}$. Hence, the products $u_{i-1}u_i$ and~$u_{i-1}'u_i'$ occurring in Theorem~\ref{5.6} are divisible by~$h^2$ in $\mathbb Z[h]$. After specializing to an even integer value of $h$, Theorem~\ref{5.6} therefore implies that~$F_{i}^a(h)$ and~$F_{i}^b(h)$ are congruent respectively to $-w_i$ and~$-w_i'$ modulo $4$; in particular each is either odd or divisible by~$4$.
\end {Remark}

\begin{Example}\label{5.7}
Recall that~$\varphi:=\tfrac{\sqrt5-1}{2}$. The special values $h=3$ and~$h=\sqrt5$ make the connection with Fibonacci-type sequences particularly transparent.
\begin{parts}
\item For~$h=3$,
\begin{equation*}\tag{Ex1}\label{eq:ex-a3}
 (a_n\psup{3})_{n\ge1}=3,8,21,55,144,377,987,\ldots,
\end{equation*}
which is $(F_4,F_6,F_8,\ldots)$ (OEIS A001906, beginning at index $2$) \cite{OEIS_A001906}, and~$\tfrac{a_{n-1}\psup{3}}{a_n\psup{3}}\to\varphi^2$ as $n\to\infty$.

\item
\begin{equation*}\tag{Ex2}\label{eq:ex-b3}
 (b_n\psup{3})_{n\ge2}=7,18,47,123,322,843,2207,5778,\ldots,
 \qquad \frac{b_{n-1}\psup{3}}{b_n\psup{3}}\to\varphi^2\qquad(n\to\infty),
\end{equation*}
The sequence $(b_n\psup{3})_{n\ge2}$ is OEIS A005248 beginning at index $2$ \cite{OEIS_A005248}.
and
\begin{equation*}\tag{Ex3}\label{eq:ex-d3}
 (d_n\psup{3})_{n\ge3}=21,54,141,369,966,2529,6621,\ldots,
 \qquad \frac{d_{n-1}\psup{3}}{d_n\psup{3}}\to\varphi^2\qquad(n\to\infty).
\end{equation*}

\item The first values of the two-index invariant are
\begin{equation*}\tag{Ex4}\label{eq:ex-Fa3}
 F_{1}^a(3),F_{2}^a(3),\ldots
 =11,80,551,3779,25904,177551,1216955,\ldots,
\end{equation*}
and
\[
 \frac{F_{i}^a(3)}{F_{i+1}^a(3)}\longrightarrow\varphi^4\qquad(i\to\infty).
\]

\item For~$h=\sqrt5$,
\begin{equation*}\tag{Ex5}\label{eq:ex-asqrt5}
 (a_n\psup{\sqrt5})_{n\ge1}
 =\sqrt5,4,3\sqrt5,11,8\sqrt5,29,21\sqrt5,76,\ldots,
 \qquad
 \frac{a_{n-1}\psup{\sqrt5}}{a_n\psup{\sqrt5}}\to\varphi\qquad(n\to\infty).
\end{equation*}
The odd and even subsequences satisfy
\[
 \frac{a_{2n-3}\psup{\sqrt5}}{a_{2n-1}\psup{\sqrt5}}\to\varphi^2\qquad(n\to\infty),
 \qquad
 \frac{a_{2n-2}\psup{\sqrt5}}{a_{2n}\psup{\sqrt5}}\to\varphi^2\qquad(n\to\infty).
\]
Moreover, Proposition~\ref{4.01} gives the exact identity
\[
 \sqrt5\,a_n\psup{3}=a_{2n+1}\psup{\sqrt5}.
\]
\end{parts}
\end{Example}

Since~$a_n\psup{2}=n+1$, the case $h=2$ represents every positive integer through this determinant sequence. The behavior changes for other values of $h$. For example, at~$h=1$, Theorem~\ref{4.4} leaves only the indecomposable virtual positive semidefinite types $A_1\psup{1}$ and~$A_2\psup{1}$, whose determinants are $1$ and~$0$, respectively; hence, the positive integer $2$ does not occur. This suggests the following problem.

\begin{Conjecture}\label{5.11}
For every positive integer $h\ne2$, there exists a positive integer $m$ such that $\det H\ne m$ for every indecomposable virtual positive semidefinite $h$-Shuhan matrix $H$.
\end{Conjecture}

The first unexplored integral case beyond the classical normalization is $h=3$.
\begin{Question}\label{5.10}
Classify the positive definite $3$-Shuhan matrices.
\end{Question}

\appendix
\section[Reproducibility check for the polynomials b-hat]{Reproducibility check for the polynomials $\widehat b_n\psup{h}$}\label{app:cas}
The determinant polynomials in Proposition~\ref{4.8} can be reproduced exactly by the following short SageMath script. This is included only as a verification aid; the proofs above use the identity $\widehat b_n\psup{h}=\det(hI_n-J_n)$ from \eqref{eq:Jn-reversal}, and do not depend on numerical experimentation.
\par\medskip
\noindent\begin{minipage}{\linewidth}
\small
\begin{verbatim}
R.<h> = PolynomialRing(QQ)

def bhat(n):
    M = matrix(R, n, n, 0)
    for i in range(n):
        M[i,i] = h
    for i in range(n-1):
        w = QQ(3)/2 if i == 0 else QQ(1)
        M[i,i+1] = -w
        M[i+1,i] = -w
    return M.det()

for n in range(2,9):
    p = bhat(n)
    print(n, p)
    print(max(r for r in p.roots(AA, multiplicities=False)))
\end{verbatim}
\end{minipage}\par\medskip

\noindent For~$n=8$, this gives
\[
 \widehat b_8\psup{h}=h^8-\frac{33}{4}h^6+\frac{85}{4}h^4
 -\frac{35}{2}h^2+\frac94,
 \qquad \mu_8\approx1.994839578989595.
\]

\section*{Acknowledgments}
Weicai Wu was supported by Guangxi Science and Technology Program under Grant No.~2026GXNSFDA\allowbreak00640036, Center for Applied Mathematics of Guangxi (Guangxi Normal University) and the Key Laboratory of Mathematical Model and Application (Guangxi Normal University).

\end{document}